\documentclass[12pt]{article}
\usepackage{amsfonts,amsmath,amsxtra}
\usepackage{latexsym}
\usepackage{amssymb}
\usepackage{color}
\usepackage[matrix,arrow,curve]{xy}

\makeatletter
\def\@seccntformat#1{\csname
the#1\endcsname\enspace} \makeatother

\def\hybrid{\topmargin 0pt      \oddsidemargin 0pt
        \headheight 0pt \headsep 0pt
        \textwidth 16.5cm
        \textheight 23cm
        \marginparwidth 0.0in
        \parskip 5pt plus 1pt   \jot = 1.5ex}
\catcode`\@=11
\def\marginnote#1{}
\newcount\hour
\newcount\minute
\newtoks\amorpm
\hour=\time\divide\hour by60 \minute=\time{\multiply\hour by60
\global\advance\minute by-\hour}
\edef\standardtime{{\ifnum\hour<12 \global\amorpm={am}%
        \else\global\amorpm={pm}\advance\hour by-12 \fi
        \ifnum\hour=0 \hour=12 \fi
      \number\hour:\ifnum\minute<10 0\fi\number\minute\the\amorpm}}
\edef\militarytime{\number\hour:\ifnum\minute<10 0\fi\number\minute}

\def\draftlabel#1{{\@bsphack\if@filesw {\let\thepage\relax
   \xdef\@gtempa{\write\@auxout{\string
      \newlabel{#1}{{\@currentlabel}{\thepage}}}}}\@gtempa
   \if@nobreak \ifvmode\nobreak\fi\fi\fi\@esphack}
        \gdef\@eqnlabel{#1}}
\def\@eqnlabel{}
\def\@vacuum{}
\def\draftmarginnote#1{\marginpar{\raggedright\scriptsize\tt#1}}

\def\draft{\oddsidemargin -0.1truein
        \def\@oddfoot{\sl RiemannZeta.tex \hfil
        \rm\thepage\hfil\sl\today\quad\militarytime}
        \let\@evenfoot\@oddfoot \overfullrule 3pt
        \let\label=\draftlabel
        \let\marginnote=\draftmarginnote
\def\@eqnnum{{\rm (\theequation)}
\rlap{\kern\marginparsep\tt\@eqnlabel}%
\global\let\@eqnlabel\@vacuum}  }

\newfont{\Bbbb}{msbm7 scaled 1\@ptsize00}
\newcommand{\zs}{\raise-1pt\hbox{$\mbox{\Bbbb Z}$}}

scaled 1\@ptsize00

\def\numberbysection{\@addtoreset{equation}{section}
        \def\theequation{\thesection.\arabic{equation}}}
\numberbysection

\renewcommand{\theequation}{\thesection.\arabic{equation}}

\def\titlepage{\@restonecolfalse\if@twocolumn\@restonecoltrue\onecolumn
     \else \newpage \fi \thispagestyle{empty}\c@page\z@
\def\thefootnote{\fnsymbol{footnote}} }
\def\endtitlepage{\if@restonecol\twocolumn \else  \fi
        \def\thefootnote{\arabic{footnote}}
        \setcounter{footnote}{0}}  
\relax \hybrid
\makeatletter
\newdimen\normalarrayskip            
\newdimen\minarrayskip               
\normalarrayskip\baselineskip \minarrayskip\jot
\newif\ifold             \oldtrue            \def\new{\oldfalse}
\def\arraymode{\ifold\relax\else\displaystyle\fi}
\def\eqnumphantom{\phantom{(\theequation)}} 
\def\@arrayskip{\ifold\baselineskip\z@\lineskip\z@
     \else
     \baselineskip\minarrayskip\lineskip1\baselineskip\fi}
\def\@arrayclassz{\ifcase \@lastchclass \@acolampacol \or
\@ampacol \or \or \or \@addamp \or
   \@acolampacol \or \@firstampfalse \@acol \fi
\edef\@preamble{\@preamble
  \ifcase \@chnum
     \hfil$\relax\arraymode\@sharp$\hfil
     \or $\relax\arraymode\@sharp$\hfil
     \or \hfil$\relax\arraymode\@sharp$\fi}}
\def\@array[#1]#2{\setbox\@arstrutbox=\hbox{\vrule
     height\arraystretch \ht\strutbox
     depth\arraystretch \dp\strutbox
width\z@}\@mkpream{#2}\edef\@preamble{\halign \noexpand\@halignto
\bgroup \tabskip\z@ \@arstrut \@preamble \tabskip\z@ \cr}%
\let\@startpbox\@@startpbox \let\@endpbox\@@endpbox
  \if #1t\vtop \else \if#1b\vbox \else \vcenter \fi\fi
  \bgroup \let\par\relax
  \let\@sharp##\let\protect\relax
  \@arrayskip\@preamble}
\def\eqnarray{\stepcounter{equation}%
              \let\@currentlabel=\theequation
              \global\@eqnswtrue
              \global\@eqcnt\z@
              \tabskip\@centering              
              \let\\=\@eqncr
              $$%
            \halign to \displaywidth  \bgroup
             \eqnumphantom \@eqnsel
      \hskip\@centering                               
    $\displaystyle  \tabskip\z@ {##}$%
    &\global\@eqcnt\@ne \hskip 2\arraycolsep
         $ \displaystyle  \arraymode{##}$\hfil
    &\global\@eqcnt\tw@ \hskip 2\arraycolsep
         $\displaystyle\tabskip\z@{##}$\hfil
         \tabskip\@centering
    &{##}\tabskip\z@\cr}
\makeatother

\def\IC{\mathbb{C}}

\def\IP{\mathbb{P}}
\def\IQ{\mathbb{Q}}
\def\IR{\mathbb{R}}
\def\IZ{\mathbb{Z}}

\def\CH {\mathcal{H}}

\def\CM {\mathcal{M}}

\def\CP {\mathcal{P}}

\def\CV {\mathcal{V}}

\def\a {{\alpha}}

\def\g {{\gamma}}
\def\s {{\sigma}}

\def\la{\lambda}

\def\e{\epsilon}

\def\wh{\widehat}

\def\Id{{\rm Id}}
\def\Fun{{\rm Fun}}

\def\c{\cdot}

\def\Mat{{\rm Mat}}

\def\ov {{\overline}}

\def\Tr{{\rm Tr}\,}

\def\<{\langle}
\def\>{\rangle}
\def\ov{\overline}

\makeatletter
\DeclareRobustCommand{\loplus}{\mathbin{\mathpalette\dog@lsemi{+}}}
\DeclareRobustCommand{\roplus}{\mathbin{\mathpalette\dog@rsemi{+}}}

\newtheorem{te}{Theorem}[section]
\newtheorem{de}{Definition}[section]
\newtheorem{prop}{Proposition}[section]           
\newtheorem{cor}{Corollary}[section]
\newtheorem{lem}{Lemma}[section]

\newcommand{\proof}{\noindent {\it Proof}. }
\newcommand\bqa{\begin{eqnarray}}
\newcommand\eqa{\end{eqnarray}}
\def\be{\begin{eqnarray}\new\begin{array}{cc}}
\def\ee{\end{array}\end{eqnarray}}
\def\beq{\begin{equation}}
\def\eeq{\end{equation}}
\def\bse{\begin{subequations}}                
\def\ese{\end{subequations}}
\def\bp{\begin{pmatrix}}
\def\ep{\end{pmatrix}}

\def\i{\imath}
\newcounter{pac}[section]

\newcounter{pacc}[subsection]

 0

\begin{document}

\title{\bf Global $GL_2$ Hecke-Baxter operator}
\author{A.A. Gerasimov, D.R. Lebedev and S.V. Oblezin}
\date{\today}
\maketitle

\renewcommand{\abstractname}{}

\begin{abstract}

  \noindent {\bf Abstract}. We construct a global Hecke-Baxter
  operator for integrable systems of arithmetic type associated with
  the group $GL_2$. It is an
  element of the global Hecke algebra associated with the double coset
  space   $GL_2(\IZ)\backslash GL_2(\IR)/O_2$.   Eigenvalues
  of the global Hecke-Baxter
  operator  acting on the $GL_2$-Eisenstein series
  are given by the   corresponding global $L$-factors. This construction generalizes
  our previous construction of the Hecke-Baxter operators  over
  local completions $\IR$ and $\IQ_p$ of the number field $\IQ$.
  Presumably, zeroes of the corresponding global $L$-factors
  should be subjected to an arithmetic version of the Bethe  ansatz
  equations.
\end{abstract}

 \vspace{5mm}


\section{Introduction}

An interpretation of a wide class of integrable systems in terms of
representation theory provides  important insights both into  theory
of integrable systems and into representation theory allowing
transferring various techniques from one  area of research into the
other. Among numerous examples we would like to mention the
formalism of the Baxter operator \cite{Ba} which was properly placed
in representation theory perspective  in \cite{GLO08} using Hecke
algebras formalism (see also \cite{G}). We coin the term Hecke-Baxter operator for a one-parameter family of elements of an appropriate Hecke algebra reproducing
  the Baxter operator for a class of  quantum integrable systems.
The  construction of the Hecke-Baxter operator for spherical
principal series representations of $GL_{\ell+1}(\IR)$ was  also
extended to general principal series representations of
$GL_{\ell+1}(\IR)$
  (for details see \cite{GLO25} and reference therein).  A remarkable fact is
that the Hecke-Baxter operators  are directly related to the
Archimedean $L$-factors attached to the corresponding
representations of $GL_{\ell+1}(\IR)$. Precisely the Archimedean
$L$-factors appear as eigenvalues of the Hecke-Baxter operators
acting on  the spherical and Whittaker functions given by specific
matrix elements of the spherical principal series  representations.
As a direct consequence,  the local $L$-factors enter integral
representations of the Whittaker functions expressed via a version
of the Gelfand-Tsetlin construction of irreducible representations
of $GL_{\ell+1}(\IR)$ \cite{GKL}. Note that the
$GL_{\ell+1}(\IR)$-Whittaker functions are eigenfunctions  of the
quantum $GL_{\ell+1}(\IR)$-Toda chains,  one of the most
well-studied finite-dimensional integrable systems associated with
the spherical principal series representations. Actually  the
Hecke-Baxter operator (depending on auxiliary parameter) provides an
alternative formulation of the quantum Toda chain. Notice in this
regard that an advantage of the Hecke algebra formulation of
integrable systems  is in a unified treatment of both continuous and
discrete  symmetries of the systems.

Not surprisingly, a proper counterpart of the Hecke-Baxter operator
exists  in the case of representation theory over non-Archimedean
fields. The case of spherical principal series representations of
$GL_{\ell+1}(\IQ_p)$ was considered in \cite{GLO08} together with
the corresponding integrable systems governing spherical and
Whittaker functions over $\IQ_p$. Connection with local $L$-factors
still holds in this case, and  the local non-Archimedean $L$-factors
show up as eigenvalues of the non-Archimedean Hecke-Baxter operators
acting on $GL_{\ell+1}(\IQ_p)$-Whittaker functions.

It is natural to expect that the Hecke-Baxter operator formalism may
be further generalized to the case of (global) number fields. This
is indeed so, and in this short note for a  compactification
$\overline{{\rm    Spec}(\IZ)}$ of ${\rm  Spec}(\IZ)$ we introduce
the corresponding  global $GL_2$ Hecke-Baxter operator acting on the
non-ramified $GL_2$-automorphic functions (functions on the double
coset $GL_2(\IZ)\backslash GL_2(\IR)/O_2$). The following result is
proven in   Theorem \ref{TH1}, Section 4. The automorphic functions
represented by specific   matrix elements of the spherical principal
series $GL_2(\IR)$-representations are eigenfunctions of the
proposed global Hecke-Baxter operator with the eigenvalues given by
the corresponding global $L$-functions generalizing the completed
form of the Riemann zeta function $\zeta(s)$:
 \be\label{zeta}
  \xi(s)\,=\,\zeta(s)\, \pi^{-\frac{s}{2}}\,\Gamma\Big(\frac{s}{2}\Big)\,.
 \ee
This result complements the constructions of \cite{JL}. In our
approach, the matrix analog of the Riemann's proof of the functional
equation for the $GL_1$ global $L$-function (given by \eqref{zeta}) naturally arises in the case of $GL_2$ (see
Appendix).

 Let us stress that in terms of quantum integrable systems we basically consider
a hyperbolic billiard on the upper half-plane modulo action of the
modular group $PSL_2(\IZ)$. This  quantum system is integrable and
is deeply connected to the quantum $GL_2(\IR)$-Toda chain: harmonics
of the quantum billiard eigenfunctions are given by solutions of the
quantum Toda chain for integer coupling constants. On the other
hand, this quantum billiard is a generalization of the Euclidean
billiard arising in the tropical limit of the $GL_2(\IR)$-Toda chain
\cite{GL}. This provides an interesting number theoretic perspective
on the tropical limit construction proposed in \cite{GL}.

One curious  point worth mentioning is as follows.
In the case of  integrable systems with discrete spectrum,
the Baxter operator is instrumental in
finding the spectrum given by common eigenvalues of quantum
Hamiltonians. Precisely,
the eigenvalues of quantum Hamiltonians are
expressed in a simple way through zeroes of the eigenvalues of
the Baxter operators (considered as functions of an auxiliary parameter). In turn,
 zeroes of the eigenvalues of the Baxter operator satisfy a set of
 equations called the Baxter equations. Our interpretation of the global
 $L$-functions as eigenvalues of the Hecke-Baxter operators points to
a possibility of existence of an  analog of the Baxter  equations in the
arithmetic setup. This might provide a new  optics
for looking at analogs of Riemann hypotheses for global $L$-functions
as well as various conjectures on
 the special values of global $L$-functions. The suggestion seems
 close to  the Faddeev-Pavlov approach \cite{FP} to studying
 analytic properties of the Riemann zeta-functions  via scattering
 theory. This line of research seems still  worth to pursue.

 Let us also note that
 the Hecke-Baxter operators (and more general elements of Hecke algebras) are
examples of averaging operators that are ubiquitous in various areas
of Mathematics and Physics. One interesting example of the averaging
operator appears in  the  Kadanov approach to the renormalization
(semi)group  in lattice quantum field theories  (see e.g.
\cite{Ka}). Fixed points of the renormalization group flow
corresponding to the eigenvalues of the Kadanov operators describe
continuum limit of the lattice theory. The analogy between the
constructions of \cite{Ka}  and of this paper is very fruitful and
will be considered in detail elsewhere. However as an obvious next
step  we  are  going to generalize the results of this note to the
global ramified  case for the groups $GL_{\ell+1}$  of higher ranks.

{\it  Acknowledgements:} The research of S.V.O. is partially
supported by the Beijing Natural Science Foundation grant IS24004.

\section{$GL_1$-automorphic forms\\ and global Hecke-Baxter operator}

In this Section we consider the almost trivial case of the Lie group $GL_1$. Our
goal is to introduce  basic  elements of the construction to proceed in the following
Sections with a more involved  case of $GL_2$.

Let us define the non-ramified $GL_1$ automorphic functions  as
functions on  $GL_1(\IR)$  invariant under the left action of
$GL_1(\IZ)$ and right action of the orthogonal subgroup $O_1\subset
GL_1(\IR)$. These functions may be considered as functions on
the double coset space
 \be\label{Mone}
  \CM_1=GL_1(\IZ)\backslash GL_1(\IR)/O_1\,,
 \ee
taking into account the following subtlety. Note that $\CM_1$  is a
$GL_1(\IZ)$-orbifold as $GL_1(\IZ)$ acts trivially on the coset
space $GL_1(\IR)/O_1$. We however are interested in the space of
functions on $\CM_1$ and thus  might ignore the orbifold structure
by  considering  functions on $GL_1(\IR)/O_1$ that are invariant
under the trivial action of $GL_1(\IZ)$. Thus taking into account
the isomorphisms
 \be
  GL_1(\IR)\simeq \IR^*, \qquad GL_1(\IZ)\simeq
  \mu_2, \qquad O_1\simeq \mu_2\,,\qquad \mu_2=\{\pm 1\}\,,
 \ee
the $GL_1$-automorphic functions may be identified with functions on
$\IR_+=\IR^*/\mu_2$.

The  double cosets space  \eqref{Mone} allows an interpretation as   a
moduli  space of circles $S^1$ supplied with  $S^1$-invariant
metrics. Indeed $\CM_1$ may be presented in the following factorized form
 \be\label{Mone1}
  \CM_1=GL_1(\IZ)\backslash GL_1(\IR)\times_{GL_1(\IR)} GL_1(\IR)/O_1\,.
 \ee
The first factor $GL_1(\IZ)\backslash GL_1(\IR)$ should be identified with
the space of lattices $L_v\subset\IR$,
 \be
  L_v=\{  nv|\,n\in \IZ\}, \qquad v\in \IR_+\,.
 \ee
By the action of $GL_1(\IR)\simeq \IR^*$ any
lattice may transformed   into the standard one, $\IZ\subset \IR$.
The second factor $GL_1(\IR)/O_1$ in \eqref{Mone1}
is identified with the  space of constant
metrics on $\IR$ with  $O_1$ being
stabilizer of a reference metric.  From this description we  infer that
the space \eqref{Mone} is  the moduli space of
 $GL_1(\IR)$-equivalence classes of  pairs of lattices  and constant metrics
on $\IR$ or equivalently as the moduli space of circles $S^1=\IR/\IZ$ supplied with
constant metrics. Algebraically $\CM_1$ may be understood as a moduli
space of rank one $\IZ$-modules $L$ supplied with a metric on its real
extension $L\otimes_\IZ \IR$.  The interpretation of $\CM_1$ as a
  moduli space of  metricized circles  provides us with a canonical
coordinate on $\CM_1$, the volume of the corresponding circle
 \be\label{VolC}
  |x|={\rm Vol}_h(\IR/L)\,.
 \ee
Here $x$ is the canonical coordinate on $GL_1(\IR)=\IR^*$ identified
with the moduli space of the oriented circles $S^1_{or}$ supplied with a constant
metric $h$. Given  a pair $(S^1_{or},h)$
we might consider corresponding volume one-form $\omega$ so that
the coordinate $x$ would be a period of this form
\be\label{VolC1}
x=\int_{\IR/L}\omega\,.
\ee

We are interested in a particular basis in the space of
$GL_1$-automorphic functions. This basis may be defined  in various
ways but having in mind  subsequent generalizations to the case of
$GL_2$ we  construct these functions via representation theory
approach. Precisely we define $GL_1$ Eisenstein functions   as
matrix elements of $GL_1(\IZ)$- and $O_1$-invariant vectors in
unitary spherical principal series
 representations of $GL_1(\IR)$.
Let $(\pi_\gamma,\CV_\gamma)$, $\gamma\in \IR$ be a  one-dimensional unitary spherical
representation of $GL_1(\IR)$, $\<\,,\,\>$ be the corresponding
 Hermitian pairing  and  $v\in \CV_{\gamma}$ be such that
$\<v,v\>=1$. By definition  unitary  spherical representations of
$GL_1(\IR)\simeq \IR^*$ are  factored through the homomorphysm $\IR^* \to
\IR_+$ and thus are  given  by
\be\label{RepGL1}
\pi_\gamma:\,\,x\longrightarrow |x|^{\imath \gamma}\,, \qquad x\in \IR^*\,.
\ee
Consider the following matrix elements  in representation $(\pi_\gamma,\CV_\gamma)$
 \be \label{ME1}
  \psi_\gamma(x)=\<v,\pi_\gamma(x)\,v\>=|x|^{\imath \gamma}\,, \qquad
  x\in GL_1(\IR)\simeq \IR^*\,.
 \ee
The function $\psi_\gamma(x)$ is  a $\mu_2$-invariant  function  on $GL_1(\IR)$ and
thus is a lift of a function on $\CM_1=\IR_+$. Corresponding function
on $\CM_1$ will be called  the $GL_1(\IR)$-Eisenstein function
associated with the representation $(\pi_\gamma,\CV_\gamma)$.
In the following we will consider
interchangeably automorphic eigenfunctions as functions on $\CM_1$
depending on $|x|$ or as $\mu_2$-invariant functions on $GL_1(\IR)$
depending on $x$.

The  $GL_1$-Eisenstein functions  may be defined also as
eigenfunctions  of appropriate operators. In the following we will
be interested in characterization of the $GL_1$-Eisenstein functions
as common eigenfunctions of elements of the Hecke algebra associated
with the space of double cosets (this formulation is especially
useful as it takes into account both the invariance under discrete
and continuous groups). Define the Hecke algebra associated with the
double coset space  \eqref{Mone} as
 a tensor product of two convolution algebras
 $\CH(GL_1(\IQ),GL_1(\IZ))$ and $\CH(GL_1(\IR),O_1)$.
Recall that Hecke algebra $\CH(G,K)$ associated with a pair
$K\subset G$ is an associative algebra of the proper subset of
$K$-biinvariant functions on $G$  under convolution. It is natural
to consider the maximal subset of functions on $G$ such that the
convolution operation is defined. In the case when $(G,K)$ is a
Gelfand pair (i.e. $K$ is a fixed subgroup of an involution of $G$)
the corresponding associative algebra is commutative. The power of the Hecke
algebra formalism is in the fact that $\CH(G,K)$ in general is not a
group algebra but replaces it in various representation theory
constructions.

The Hecke algebra $\CH(GL_1(\IR),O_1)$, the  algebra of
$O_1$-biinvraint functions on $GL_1(\IR)$,  acts naturally on the
functions on $GL_1(\IR)/O_1$   and in particular on the functions on
double coset $\CM_1$ via convolution.   Note that it  does not take
into account the $\IZ$-structure responsible for the lattice moduli
space interpretation of $\CM_1$. To take into account this
arithmetic structure we consider another Hecke algebra
$\CH(GL_1(\IQ),GL_1(\IZ))$ which we identify with the convolution
algebra of $GL_1(\IZ)$-biinvariant generalized functions on
$GL_1(\IR)$ supported at $GL_1(\IQ)\subset GL_1(\IR)$. It is easy to
verify that the algebras $\CH(GL_1(\IQ),GL_1(\IZ))$ and
$\CH(GL_1(\IR),O_1)$
 are (mutually) commutative associative algebras acting on
functions on $\CM_1$ from the left and from the right
correspondingly. Note that in the considered case of $GL_1$ the Hecke algebras are
actually the group algebras of the quotient groups $GL_1(\IQ)/GL_1(\IZ)$ and
$GL_1(\IR)/O_1$.

It is instructive to describe the $\CH(GL_1(\IQ),GL_1(\IZ))$-action
 considering $\CM_1$ as a  moduli space of
 metricized circles. We define  operations $T_{p/q}$, $p/q\in \IQ_+^*$ on lattices as follows.
For a given lattice $L$, we first  take a lattice $L_{1/q}$ such
that $L\subset L_{1/q}$ and $[L_{1/q}:L]=q$. Then we consider a
sublattice $L_{p/q}\subset L_{1/q}$ of index $[L_{1/q}:L_{p/q}]=p$.
Combining these operations we define the following operator acting
on functions on the space of lattices:
  \be
   (T_{p/q}\cdot f)(L)\,
   =\sum_{L\subset L_{1/q}\supset L_{p/q} }\,f(L_{p/q}),
   \qquad
   [L_{1/q}:L]=q, \quad [L_{1/q}:L_{p/q}]=p\,.
  \ee
 In terms of functions of $|x|\in \IR_+$ this
reduces to a simple mupltiplication operation
 \be
  (T_{p/q}\cdot f)(|x|)=f(p|x|/q)\,.
 \ee
These operators belong to the Hecke algebra
$\CH(GL_1(\IQ),GL_1(\IZ))$ and satisfy the following relations:
 \be
  T_{p_1/q_1}\circ T_{p_2/q_2}=T_{(p_1p_2)/(q_1q_2)}\,.
 \ee
Let us remark that the collection of operators $T_{p/q},\,p/q\in\IQ^*_+$
provides a $GL_1$-analog of the modular tower structure arising in the case of
$GL_2$.

To construct  a meaningful generating function   we consider a
multiplicative semigroup ${\rm Mat}^*_1(\IZ) \subset GL_1(\IQ)$
 of integer non-zero one by one matrices
\be
{\rm Mat}^*_1(\IZ)={\rm Mat}_1(\IZ)\cap
GL_1(\IQ)\,,\qquad {\rm Mat}^*_1(\IZ)=\IZ-\{0\}\,.
\ee
The semigroup is acted by $GL_1(\IZ)$ and the quotient may be
identified with  $\IZ_+$:
\be\label{Subs}
GL_1(\IZ)\backslash {\rm Mat}_1^*(\IZ)=\IZ_+\,.
\ee
Now we consider
elements of the Hecke algebra  $\CH(GL_1(\IQ),GL_1(\IZ))$
supported at the subset \eqref{Subs} and acting
via
 \be\label{Tn}
  (T_n\cdot f)(L)=f(nL)\,,\qquad n\in \IZ_+\,.
 \ee
Equivalently, in terms of functions on $GL_1(\IR)$ we have
 \be\label{Tnx}
  (T_n\cdot f)(x)=f(nx)\,,\qquad n\in \IZ_+\,, \qquad x\in \IR^*\,.
 \ee
These operators may be conveniently combined into  the generating series
 \be\label{GenZ}
  Q^{GL_1(\IZ)}_s\,=\sum_{n=1}^{\infty}\frac{1}{n^s} T_n\,,
 \ee
with its action on the functions of $x$  given by
 \be
  (Q^{GL_1(\IZ)}_s\cdot f)(x)\,=\sum_{n=1}^{\infty} n^{-s}\,f(nx)\,.
  \ee
This may be written in the following form  allowing a direct  generalization
to the case of $GL_2$ in the next Section
 \be
  (Q^{GL_1(\IZ)}_s\cdot f)(x)\,=\sum_{n\in GL_1(\IZ)\backslash {\rm Mat}_1^*(\IZ)}
  n^{-s}\,f(nx)\,.
 \ee

Introduce a kind of generating function for elements of the Hecke
algebra $\CH(GL_1(\IR),O_1)$ providing a proper counterpart for the
generating function \eqref{GenZ}. Such generating functions (for
more general case of $GL_{\ell+1}(\IR)$) were first introduced in
\cite{GLO08}  under the name of the Hecke-Baxter operator. Precisely
the $GL_1(\IR)$ Hecke-Baxter operator is the integral operator,
 \be
  (Q^{GL_1(\IR)}_s\cdot f)(x)\,
  =\int\limits_{\IR^*}\!d\mu^G_{\IR^*}(y)\,\, |y|^s\,f(y^{-1}x)\,,\quad
  d\mu^G_{\IR^*}(y)=e^{-\pi y^2}\,\frac{dy}{y}\,,
 \ee
acting by convolution with the following $O_1$-biinvariant function
on $GL_1(\IR)$
 \be\label{ArKernel}
  Q^{GL_1(\IR)}_s(y)=|y|^s\,e^{-\pi y^2}\,.
 \ee
In the following,  for brevity, we identify suitable functions on
Lie groups, operators
 obtained by the actions of these functions via convolution and
 the corresponding integral kernels.

\begin{prop} The matrix elements  \eqref{ME1}
 \be\label{ME11}
  \psi_{\gamma}(x)=\<v,\pi_{\gamma}(x)\,v\>=|x|^{\imath \gamma}\,,
 \ee
are common  eigen-functions of the operators $Q^{GL_1(\IZ)}_s$ and\,
 $Q^{GL_1(\IR)}_s$:
 \be
  \bigl(Q^{GL_1(\IZ)}_s\cdot  \psi_{\gamma}\bigr)(x)\,
  =\,\zeta(s-\imath \gamma)\, \psi_{\gamma}(x)\,,
 \ee
 \be
  \bigl(Q^{GL_1(\IR)}_s\cdot  \psi_{\gamma}\bigr(x)\,
  =\,L^{\IR}(s-\imath \gamma)\,\psi_{\gamma}(x)\,.
 \ee
The  eigenvalues are given by
 \be
  \zeta(s)=\sum_{n\in \IZ_+}\,\frac{1}{n^s}\,,\qquad
  L^{\IR}(s)=\pi^{-\frac{s}{2}}\,\Gamma\left(\frac{s}{2}\right)\,,
 \ee
where we impose ${\rm Re}(s)>1$ for convergence.
\end{prop}

\proof Using
 \be
  (T_n\cdot \psi_{\gamma})(x)\,
  =\,\psi_{\gamma}(nx)\,
  =\,n^{\imath \gamma}\,
  \psi_{\gamma}(x)\,,
 \ee
we indeed find out
 \be
  \bigl(Q^{GL_1(\IZ)}_s\cdot  \psi_{\gamma}\bigr)(x)\,
  =\Big(\sum_{n\in\IZ_+}\frac{1}{n^{s-\imath \gamma}}\Big)\,
  \psi_{\gamma}(x)\,
  =\,\zeta(s-\imath \gamma)\, \psi_{\gamma}(x)\,.
 \ee
For the  Archimedean Hecke-Baxter operator $Q^{GL_1(\IR)}_s$ acting
on the matrix element \eqref{ME1} via
 \be
  \bigl(Q^{GL_1(\IR)}_s\cdot \psi_{\gamma}\bigr)(x)\,
  =\int\limits_{\IR^*}\!\frac{dy}{y}\,\,|y|^s\,e^{-\pi y^2}\,
  \psi_{\gamma}(y^{-1}x)\,
  =\,\pi^{-\frac{s-\imath \gamma}{2}}\,
  \Gamma\left(\frac{s-\imath \gamma}{2}\right)\,\psi_{\gamma}(x)\,,
 \ee
the analogous statement basically reduces to the integral
representation of the Gamma-function.  $\Box$

Now we introduce  the  main object of our considerations in this Section,
global Hecke-Baxter operator $\wh{Q}_s$.

\begin{de} The $GL_1$ global Hecke-Baxter operator is the operator
  acting in the space of functions on
  $GL_1(\IR)\simeq \IR^*$ via convolution with the following  function
 \be\label{GLB}
  \wh{Q}^{GL_1}_s(x)\,
  =\,\frac{1}{2}\,|x|^s\Big(\Theta(0|\imath
  x^2)-1\Big)\,,\qquad x\in \IR^*\,,
\ee
where the theta-constant is given by
 \be\label{Tconst}
  \Theta(0|\tau)=\sum_{n\in \IZ} e^{\imath \pi \tau n^2}\,.
  \ee
\end{de}

Along with the theta-constant
\be\label{Tconst1}
  \Theta(0|\tau)=\sum_{n\in {\rm Mat}_1(\IZ)} e^{\imath \pi \tau n^2}\,,
  \ee
it is useful to introduce the following modified theta-series
 \be\label{Tconst12}
  \Theta^*(0|\tau)\,
  =\!\!\sum_{n\in {\rm Mat}^*_1(\IZ)}\!\!
  e^{\imath \pi \tau n^2}\,
  =\, \Theta(0|\tau)-1\,,
 \ee
 \be\label{Tconst123}
  \Theta^{**}(0|\tau)\,
  =\!\!\sum_{n\in GL_1(\IZ)\backslash {\rm Mat}^*_1(\IZ)}\!\!
  e^{\imath \pi \tau n^2}\,
  =\,\frac{1}{|GL_1(\IZ)|}\Big(\Theta(0|\tau)-1\Big)\,.
  \ee

\begin{prop}  Consider the  matrix element \eqref{ME1},
 \be \label{ME12}
  \psi_{\gamma}(x)=\<v,\pi_{\gamma}(x)\,v\>=|x|^{\imath \gamma}\,,
 \ee
in the unitary spherical principal series representation
$(\pi_\gamma,\CV_\gamma)$ of
$GL_1(\IR)$. Define the completed zeta-function by
 \be\label{ExtRM}
  \xi(s)\,=\,\zeta(s)\,\pi^{-\frac{s}{2}}\,\Gamma\Big(\frac{s}{2}\Big)\,.
 \ee
Then the global Hecke-Baxter operator \eqref{GLB} acts on \eqref{ME12}
via multiplication by a shifted completed zeta-function
 \be\label{GL1act}
  \bigl(\wh{Q}^{GL_1}_s\cdot \psi_\gamma\bigr)(x)\,
  =\xi(s-\imath \gamma)\,\psi_\gamma(x)\,,\qquad{\rm Re}(s)>1\,.
 \ee
\end{prop}

\proof We have
 \be
  \bigl(\wh{Q}^{GL_1}_s*\psi_\gamma\bigr)(x)\,
  =\!\int\limits_{\IR^*}\frac{dy}{y}\,\, |y|^{s}\,\,\,\Theta^*(0|\imath  y^2)\,\psi_\gamma(y^{-1}x)\\
  =\sum_{n\in \IZ_+}\,\int\limits_{\IR^*}\!\frac{dy}{y}\,|y|^{s}\,e^{-\pi |ny|^2}\,
  |y|^{-\imath \gamma}\,\psi_\gamma(x)\,\\
  =\,\pi^{-\frac{s-\imath \gamma}{2}}\,\Gamma\Big(\frac{s-\imath \gamma}{2}\Big)\,
  \left(\sum_{n\in \IZ_+} \frac{1}{n^{s-\imath
        \gamma}}\right)\,\,\psi_\gamma(x)\\
        =\,\pi^{-\frac{s-\imath \gamma}{2}}\,\Gamma\Big(\frac{s-\imath \gamma}{2}\Big)\,
  \zeta(s-\imath \gamma)\,\,\psi_\gamma(x)\,,
\ee
thus arriving at the required identity. $\Box$

Let us notice that in the simple case of the trivial representation
$(\pi_{\g=0},\CV_{\gamma=0})$ the identity
\eqref{GL1act} reduces to the standard
integral expression for the completed zeta-function:
 \be\label{ZetaInt}
  \xi(s)\,
  =\,\int\limits_{\IR_+}\frac{dt}{t}\,\,
  t^{\frac{s}{2}}\,\,\sum_{n\in \IZ_+}\,e^{-\pi tn^2}
  =\,\int\limits_{GL_1(\IR)/GL_1(\IZ)}\frac{dt}{t}\,\,
  t^{\frac{s}{2}}\,\,\Theta^*(0|\imath t) \,.
 \ee

The fundamental property of the completed  Riemann zeta-function
\eqref{ExtRM} is the functional relation
 \be\label{FR1}
  \xi(1-s)=\xi(s)\,.
 \ee
Its proof is standard and goes back to Riemann. First we
 decompose the integral
 \be
  \xi(s)\,
  =\,\int\limits_0^{\infty}\frac{dt}{t}\,\,
  t^{\frac{s}{2}}\,\,\Theta^*(0|\imath t)\,
  =\int\limits_0^{\infty}\frac{dt}{t}\,\,
  t^{\frac{s}{2}}\,\Big(\Theta(0|\imath t)-1\Big)\,,
 \ee
as follows
 \be
  \xi(s)\,=\,-\frac{1}{2}\int\limits_0^1\frac{dt}{t} \,t^{\frac{s}{2}}\,
  +\,\frac{1}{2}\int_0^{1}\frac{dt}{t}\,t^{\frac{s}{2}}\,\,\Theta(0|\imath t)\\
  +\,\frac{1}{2}\int\limits_1^{\infty}\frac{dt}{t}\,t^{\frac{s}{2}}\,\,\Big(\Theta(0|\imath t)-1\Big)\,.
 \ee
Applying the modular transformation properties of the theta-constant
 \be\label{ModTr}
  \Theta(0|-\tau^{-1})=(-\imath \tau)^{1/2}\,\Theta\left(0|\tau\right)\,,
 \ee
for the second term we obtain
 \be
  \xi(s)\,=-\frac{1}{2}\int\limits_0^1\frac{dt}{t}\,t^{\frac{s}{2}}\,
  +\, \frac{1}{2}\int_1^{\infty}\frac{dt}{t}\,t^{\frac{1-s}{2}}\,\,\Theta(0|\imath t)\\
  +\,\frac{1}{2}\int\limits_1^{\infty}\frac{dt}{t}\,t^{\frac{s}{2}}\,\,\Big(\Theta(0|\imath t)-1\Big)\\
  =-\left(\frac{1}{s}\,+\,\frac{1}{1-s}\right)\,
  +\,\frac{1}{2}\int\limits_1^{\infty}\frac{dt}{t}\,(t^{\frac{1-s}{2}}+t^{\frac{s}{2}})\,\,
  \Big(\Theta(0|\imath t)-1\Big)\,.
 \ee
This presentation is explicitly invariant under the inversion $s\mapsto
1-s$ and allows an analytic continuation over $s$, therefore it
verifies the functional equation.

The  matrix elements \eqref{ME12} also respect an appropriate reflection symmetry
 \be
  \psi_{-\gamma}(x^{\tau})=\psi_\gamma(x)\,,
 \ee
where $x^{\tau}:=x^{-1}$ is the involution on the group
$GL_1(\IR)$. Taking into account that matrix elements \eqref{ME12} are
eigenfunctions of the Hecke-Baxter operator $\hat{Q}_s^{GL_1}$ with
the eigenvalues expressed through completed zeta-function  one
expects that the kernel of the Hecke-Baxter integral operator
should also  satisfy a form of functional equation. Indeed we have the
following relation
\be\label{FR2}
\wh{Q}^{GL_1}_{1-s}(x^\tau)\,+\,\frac{1}{2}|x^\tau|^{1-s}\,
=\,\wh{Q}^{GL_1}_{s}(x)\,+\,\frac{1}{2}|x|^{s}\,,
\ee
where the terms $|x|^s$ compensate the correction terms
entering the expression \eqref{GLB} of the kernel via theta-constant.
The functional relation  \eqref{FR2} is a direct consequence of the
modular properties \eqref{ModTr} of the theta-constant \eqref{Tconst}.
Thus we have a deep connection between properties
of the global Hecke-Baxter  operator and analytic properties of the
completed  Riemann zeta-function.

It is possible to interpolate between the global and Archimedean
Hecke-Baxter operators via considering a $GL_1$-analog of the
congruence (semi)groups. For $N\geq0$, let us introduce the
following semigroup $\IZ_+^{(N)}\subset\IQ^*_+$:
 \be
  \IZ_+^{(N)}\,=\,\bigl\{\eta \in \IQ_+^*\,\bigr|\,\,\eta=1+Nm, \,m\in \IZ_{\geq 0}\bigr\}\,,\quad
  \IZ_+^{(0)}=\{1\}\,.
 \ee
Then consider the following generating function of elements of the
Hecke algebra \\$\CH(GL_1(\IQ),GL_1(\IZ))$:
 \be
  Q^{GL_1(\IZ)}_{s,N}=\sum_{n=0}^{\infty} \frac{1}{(1+nN)^s}\,T_{1+nN}\,.
 \ee
Therefore the  modified kernel of the global Hecke-Baxter operator is
given by
 \be\label{GLBN}
  \wh{Q}^{GL_1}_{s,N}(x)\,
  =\,\frac{1}{2}\,|x|^s\,\,\Theta^{(N)}(0|\imath x^2)\,,\qquad x\in
  \IR^*\,,\qquad N>1\,.
 \ee
Here
\be
\Theta^{(N)}(0|\tau)=\sum_{n\in \IZ} e^{\imath \pi \tau (1+Nn)^2}=
e^{\imath \pi \tau }\,
\sum_{n\in \IZ} e^{\imath \pi \tau N^2n^2+2\pi \imath Nn\tau}=
e^{\imath \pi \tau }\,\Theta_{N^2}\Big(\frac{\tau}{N}\Big|\tau\Big)\,,
\ee
where the level $k$ theta function is defined by:
 \be
  \Theta_k(z|\tau)\,
  =\sum_{n\in \IZ} e^{\imath \pi k\tau n^2+2\pi \imath knz}\,.
 \ee
Now it is easy to check that by taking the limit $N\to+\infty$, the kernel \eqref{GLBN}
of the modified global Hecke-Baxter operator   turns into the kernel
\eqref{ArKernel} of the Archimedean Hecke-Baxter operator.  This
provides a kind of regularization of the Archimedean  Hecke-Baxter
operator.

\section{$GL_2$-automorphic forms }

In this Section we recall a  construction of the Eisenstein
functions for $GL_2$ (for a review see e.g. \cite{ILP} and references therein). 
Let us start with considering the  double coset space
 \be\label{Mtwo}
  \CM_2=GL_2(\IZ)\backslash GL_2(\IR)/O_2\,.
 \ee
The space \eqref{Mtwo}  is an orbifold and we will define the space
of functions on \eqref{Mtwo} as functions on $GL_2(\IR)/O_2$
invariant under the $GL_2(\IZ)$-action from the left. The double
coset space  $\CM_2$ allows an interpretation as a moduli space of
real two-tori $T^2$ supplied with $T^2$-invariant metrics. Indeed,
the space $\CM_2$ may be identified with the space of pairs of
lattices $L\subset\IR^2$ with constant metrics $h$ modulo the
simultaneous action of $GL_2(\IR)$. Taking into account that via
linear transformations any lattice in $\IR^2$ may be transformed to
the standard one $\IZ^2\subset \IR^2$ we arrive at the
identification of $\CM_2$ with the moduli space of metricizes
two-tori. Constant metric on $T^2$ defines a conformal structure,
and therefore a complex structure, supplying $T^2$ with a structure
of elliptic curve $E(\IC)$. As a result the space $\CM_2$ is
naturally fibred over the moduli space $\CM_2^c$ of elliptic curves.
A fiber of the projection $\CM_2\to \CM_2^c$ may be identified with
$\IR_+$ supplied with the natural coordinate, the volume of
$T^2=\IR^2/\IZ^2$ in the considered metric. The  space $\CM_2^c$ of
complex structures  has the double coset description as the upper
complex half-plane
 \be\label{ISO12}
  \CH_+=PSL_2(\IR)/SO(2)\,,
 \ee
modulo action of the discrete group $PSL_2(\IZ)$ from the left.
Using the standard linear coordinates on $\CH_+=\{\tau\in \IC|{\rm
Im}(\tau)>0\}$ the isomorphism \eqref{ISO12} may be described by the
following identification of $\CH_+$ with the space of $SO_2$-cosets:
 \be
  \tau=(\tau_1+\imath \tau_2)\in\CH_+\, \longmapsto\,
  \frac{1}{\sqrt{\tau_2}}
  \begin{pmatrix}  \tau_2 & \tau_1\\ 0 & 1\end{pmatrix}SO_2\,
  \subset\,PSL_2(\IR)\,,\qquad  \sqrt{\tau_2}>0\,.
 \ee
The $PSL_2(\IZ)$-action from the left on $\CH_+$ is realized by the
fraction-linear transformations.

Now we introduce a special kind of $GL_2$-automorphic functions, the
$GL_2$-Eisenstein functions. The $GL_2$-Eisenstein functions are
associated with  spherical principal series representations entering
the (continuous spectrum part of) decomposition of the
$GL_2(\IR)$-representation in the space of functions on
$GL_2(\IZ)\backslash GL_2(\IR)$ with the action of $GL_2(\IR)$ from
the right:
 \be
  \bigl(\pi(g)\cdot f\bigr)(\tilde{g})=f(\tilde{g}\cdot g),\qquad f\in
  {\rm Fun}(GL_2(\IZ)\backslash GL_2(\IR))\,.
  \ee
The irreducible components corresponding to spherical principal
series representations are in one to one correspondence with the
elements of ${\rm Fun}(GL_2(\IZ)\backslash GL_2(\IR))$ invariant
under the action of the subgroup $O_2\subset GL_2(\IR)$. This
correspondence follows from  the uniqueness of spherical vectors in
spherical principal series representations (see e.g. \cite{GGPS}).
As is usual for continuous part of spectral decomposition, the
matrix elements of irreducible constituents do not belong to the
space of the square integrable functions $L^2(GL_2(\IZ)\backslash
GL_2(\IR))$ but should be understood as half-densities on the
product of the Lie group and its unitary dual. Our considerations
will be local over the unitary dual space, and in the following we
will ignore this subtlety considering matrix elements as elements of
$L^2(GL_2(\IZ)\backslash GL_2(\IR))$.

 For $\g=(\g_1,\g_2)\in\IR^2$, let $(\pi_\gamma,\CV_\gamma)$ be a unitary spherical
 principal series  representation  of $GL_2(\IR)$
realized via induction
 \be
  \pi_{\g}={\rm Ind}_B^{GL_2(\IR)}\chi^B_\gamma\,,
\ee
from the Borel subgroup $B\subset GL_2(\IR)$
(identified with the subgroup of  lower triangular matrices) via the
spherical character of $B$
\be
  \chi^B_\gamma(b)=\prod_{j=1}^2 |b_{jj}|^{\imath \g_j-\rho_j}\,,\quad
  \rho=\Big(\frac{1}{2},\,-\frac{1}{2}\Big)\,.
 \ee
The representation space $\CV_\gamma$
 \be
  \CV_{\g}=\bigl\{f\in\Fun(GL_2(\IR))\,\bigr|\quad
  f(bg)=\chi^B_{\g}(b)\,f(g),\,\,b\in B\bigr\}\,,
 \ee
 supports the $GL_2(\IR)$-action from the right. The
 representation $(\pi_\gamma,\CV_\gamma)$ may be realized in the
space of functions on $B\backslash GL_2(\IR)=\IP^1(\IR)$, which in
turn can be identified with the compactification of the (opposite)
unipotent subgroup $N_+\subset GL_2(\IR)$:
 \be
  N_+\,=\,\Big\{n_x=\Big(
  \begin{smallmatrix} 1&&x\\&&\\0 &&1\end{smallmatrix}\Big)\,\Big|
                                                 \quad x\in\IR\Big\}\,.
 \ee
Explicitly, the $GL_2(\IR)$ action in $\CV_\gamma\subset
L^2(B\backslash GL_2(\IR))$ is given by, for
$g=\Big(\begin{smallmatrix}a&&b\\&&\\c&&d\end{smallmatrix}\Big)$,
 \be\label{TrRules}
  [\pi_\gamma(g)\, f](x)\,=\,f(n_xg)\,
 =\,|\det g|^{\imath\gamma_2+\frac{1}{2}}\,|a+xc|^{\imath (\gamma_1-\gamma_2)-1}\,
  f(g\cdot x)\,,
 \ee
providing the following $GL_2(\IR)$-action on $\IP^1(\IR)$:
 \be\label{RLaction}
  g\cdot x\,=\,\frac{b+xd}{a+xc}\,,\qquad
  g=\Big(\begin{smallmatrix}a&&b\\&&\\c&&d\end{smallmatrix}\Big)\,.
  \ee

The Hilbert space structure on $\CV_{\g}$ is defined via the
pairing,
 \be\label{PAIR}
  \<\phi_1\,,\phi_2\>\,=\int\limits_{\IR}\!dx\,\,\,\ov{\phi_1(x)}\,\,\phi_2(x)\,.
  \ee
We would like to represent the $GL_2$-Eisenstein functions in terms
of matrix elements of the spherical principal series representation
$(\pi_{\g},\CV_{\g})$. As we will see the corresponding matrix
elements are not well-defined for the unitary principal series and
require analytic continuation of the representation parameters
$\gamma=(\gamma_1,\gamma_2)\in \IR^2$. In turn, this implies a
replacement of the structure of Hilbert space $\CV_\gamma$ by a pair
of a space and its  dual. Precisely, we  supply the Hilbert space
$\CV_{\gamma}$ with a structure of the rigged Hilbert spaces
$\CV_{\gamma}^{(t)}\subset \CV_{\gamma}\subset \CV_{\gamma}^{(g)}$
(the Gelfand triple), where $\CV_{\gamma}^{(t)}$ is the subspace of
smooth test functions and $\CV_{\gamma}^{(g)}$ is the  space of
generalized functions (tempered distributions). The pairing
\eqref{PAIR} on $\CV_{\g}$ extends to the duality pairing between
$\CV_{\gamma}^{(t)}$ and $\CV_{\gamma}^{(g)}$. In the following we consider
$\CV_{\gamma}^{(t)}$ to be the  $GL_2(\IR)$-module with
representation parameters $\gamma=(\gamma_1,\gamma_2)\in \IC^2$ having non-zero  
imaginary parts. The dual space $\CV_{\bar{\g}}^{(g)}$ then has the
structure of a $GL_2(\IR)$-module
with the complex conjugated  representation parameters $\bar{\gamma}=(\bar{\g}_1,\bar{\g}_2)$.

To construct a matrix element representation of the
$GL_2$-Eisenstein functions we start with explicit construction of
$GL_2(\IZ)$- and $O_2$-invariant vectors. Let
$\phi_{O_2}\in\CV_{\g}^{(t)}$ be a spherical vector (i.e. vector
invariant under the action of $O_2\subset GL_2(\IR)$), and let
$\phi^{\vee}_{GL_2(\IZ)}$ be a $GL_2(\IZ)$-invariant vector in
$\CV_{\bar{\g}}^{(g)}$.

\begin{lem}
In the representation $(\pi_\gamma,\CV_\gamma)$ given by \eqref{TrRules}
  the $O_2$-invariant vector $\phi_{O_2}\in\CV_{\g}^{(t)}$ and
the  $GL_2(\IZ)$-invariant vector
$\phi^{\vee}_{GL_2(\IZ)}\in\CV_{\bar{\g}}^{(g)}$ are unique (up to normalization)
and may be chosen in the following form:
 \be\label{SphVect1}
  \phi_{O_2}(x)\,
  =\,(1+x^2)^{\frac{\imath(\gamma_1-\gamma_2)-1}{2}}\,,
 \ee
 \be \label{SphVect2}
  \phi^{\vee}_{GL_2(\IZ)}(x)=\sum_{(m,n)\in \CP}\!\!
  |n|^{-\imath(\bar{\gamma}_1-\bar{\gamma}_2)}\, \delta(m+nx)\,,
  \ee
where $\delta(x)$ is the delta-function, and
  \be\label{Range}
  \CP\,=\,\bigl\{ (m,n)\in \IZ^2\setminus\{(0,0)\}\,\bigl|\quad
  {\rm gcd}(m,n)= 1\bigr\}/\sim \,,
  \ee
 with  the following equivalence relation
\be
(m,n)\sim (-m,-n)\,.
 \ee
 \end{lem}

\proof Elements of $O_2\subset GL_2(\IR)$
may be written in the following form
 \be
  k(\theta)=\begin{pmatrix} \cos \theta & (-1)^\epsilon \sin \theta
    \\-\sin\theta &(-1)^\epsilon\,   \cos\theta
  \end{pmatrix} \in O_2\,,\quad0\leq\theta<2\pi\,, \quad \epsilon\in \{0,1\}\,.
 \ee
By \eqref{RLaction}, a direct calculation gives
  \be\label{Tr123}
  1+\bigl(k(\theta)\cdot x\bigr)^2\,
  =\,1+\Big(\frac{(-1)^{\e}\sin\theta+x(-1)^{\e}\cos\theta}
  {\cos\theta-x\sin\theta}\Big)^2\,
  =\,\frac{1+x^2}{(\cos\theta-x\sin\theta)^2}\,.
  \ee
Taking into account $|\det k(\theta)|=1$ we infer from \eqref{TrRules} the
$O_2$-invariance of the vector \eqref{SphVect1}. Uniqueness of the
spherical vector follows from  the fact that it should depend on
$x^2$ (due to invariance  under the action of diagonal elements of $O_2$) and the
transformation properties \eqref{Tr123}.

Next,   we find a $GL_2(\IZ)$-invariant vector
$\phi^{\vee}_{GL_2(\IZ)}\in\CV_{\bar{\g}}^{(g)}$ by solving the   following
equation,
 \be
  [\pi^{\vee}_\gamma(g)\,\phi^{\vee}_{GL_2(\IZ)}](x)\,
  =\,\phi^{\vee}_{GL_2(\IZ)}(x)\,,\qquad
  g\in GL_2(\IZ)\,.
 \ee
 The group $GL_2(\IZ)$ is generated by
 \be
  T=\Big(\begin{smallmatrix}1&&1\\&&\\0&&1\end{smallmatrix}\Big),\qquad
  S=\Big(\begin{smallmatrix}0&&1\\&&\\1&&0\end{smallmatrix}\Big),\qquad
  R=\Big(\begin{smallmatrix}-1&&0\\&&\\0&&1\end{smallmatrix}\Big)\,.
 \ee
By \eqref{RLaction} the generators are acting via
 \be\label{GL2Zrel}
  [\pi^{\vee}_{\g}(T)\,f](x)\,=\,f(x+1),\qquad
  [\pi^{\vee}_{\g}(R)\,f](x)\,=\,f(-x),\\
  \,[\pi^{\vee}_{\g}(S)\,f](x)\,=\,|x|^{\i(\bar{\g}_1-\bar{\g}_2)-1}\,f(x^{-1})\,.
 \ee
Considering the expression \eqref{SphVect2}, its $R$-invariance
reduces to the change $n\mapsto -n$ of the summation variable,  and
invariance under the $T$-action can be verified by the  change of the variable
$m\mapsto m-n$ (which does not spoil the condition ${\rm gcd}(m,n)=1$).
To check the invariance with respect to $S$-action we take into account the following
identity
 \be
  |x|^{\imath (\bar{\g}_1-\bar{\g}_2)-1}\, \delta(m+x^{-1}n)\,
  =\,|x|^{\imath (\bar{\g}_1-\bar{\g}_2)}\, \delta(xm+n)\,
  =\,\frac{|n|^{\imath (\bar{\g}_1-\bar{\g}_2)}}
  {|m|^{\imath (\bar{\g}_1-\bar{\g}_2)}}\, \delta(xm+n)\,.
 \ee
 This completes a verification of the
required properties of \eqref{SphVect1} and \eqref{SphVect2}. $\Box$

The summation set $\CP$ in \eqref{Range} allows for an
interpretation as a set of cosets  of the group $GL_2(\IZ)$  with
respect to an appropriate subgroup. To motivate this let us write
down the formal expression for $GL_2(\IZ)$-invariant
 vector $\phi$ obtained by averaging over the $GL_2(\IZ)$-action of some vector
$\phi^{(0)}\in \CV_{\bar{\g}}^{(g)}$:
 \be\label{GL2Zvec}
  \phi(x)\,
  =\sum_{\alpha\in GL_2(\IZ)/{\rm St}_{\phi^{(0)}}}\!\!
  [\pi^{\vee}_\gamma(\alpha)\,\phi^{(0)}](x)\\
  =\sum_{\alpha\in GL_2(\IZ)/{\rm St}_{\phi^{(0)}}}\!\!
  \bigl|k+xl\bigr|^{\imath (\bar{\g}_1-\bar{\g}_2)-1}\,
  \phi^{(0)}\left(\frac{m+ xn}{k+xl}\right)\,,
 \ee
where
 \be
  \alpha\,=\,\Big(\begin{smallmatrix} 
  k&& m\\&&\\l&& n\end{smallmatrix}\Big)\in
  GL_2(\IZ), \qquad |\det g|=1\,,
 \ee
and ${\rm St}_{\phi^{(0)}}$ is the stabilizer of $\phi^{(0)}$ in
$GL_2(\IZ)$. Note that here we have the right cosets
  $GL_2(\IZ)/{\rm St}_{\phi^{(0)}}$ as
 \be
  [\pi^{\vee}_\gamma(\alpha\,\alpha')\,\phi^{(0)}](x)\,
  =\,[\pi^{\vee}_\gamma(\alpha)\,
  \phi^{(0)}](x)\,\qquad \alpha'\in {\rm
  St}_{\phi^{(0)}}\subset GL_2(\IZ)\,.
 \ee
Let us choose as an   initial function the delta-function
 \be
  \phi^{(0)}(x)=\delta(x)\,,
 \ee
which is invariant under the action of the subgroup
 \be
  B(\IZ)=\Big\{\beta=\Big(\begin{smallmatrix}
  (-1)^{\epsilon_1} && 0 \\&&\\
  r && (-1)^{\epsilon_2}\end{smallmatrix}\Big)\,\Big|\,\,
  \epsilon_{1,2}\in\{0,1\}, r\in
  \IZ\Big\}\,\subset
GL_2(\IZ)\,.
 \ee
Indeed, by the explicit expression of the $GL_2(\IZ)$-action, for
$\beta=\Big(\begin{smallmatrix}
  (-1)^{\epsilon_1} && 0 \\&&\\
  r && (-1)^{\epsilon_2}\end{smallmatrix}\Big)$,
 \be
  [\pi^{\vee}_\gamma(\beta)\,\delta](x)=
  \Big|(-1)^{\e_1}+xr\Big|^{\imath (\gamma_1-\gamma_2)-1}\,
  \delta\left(\frac{(-1)^{\e_2}x}{(-1)^{\e_1}+xr}\right)\,
  =\,\delta(x)\,.
 \ee
Thus the  sum in the formal expression shall be reduced to summation
over the space of cosets $GL_2(\IZ)/B(\IZ)$. As a result, for a
$GL_2(\IZ)$-invariant vector \eqref{GL2Zvec} we obtain 
 \be
  \phi^{\vee}_{GL_2(\IZ)}(x)\,
  =\sum_{\alpha\in GL_2(\IZ)/B(\IZ)}\,
  \,[\pi^{\vee}_\gamma(\alpha)\,\delta](x)\,.
 \ee
To establish an equivalence of this  expression with  \eqref{SphVect2}
we use the following explicit  description
of the coset representatives of $GL_2(\IZ)/B(\IZ)$.

\begin{lem}\label{BZG}
For the subgroup $B(\IZ)=B\cap GL_2(\IZ)$,
 \be
  B(\IZ)=\Big\{\beta=\Big(\begin{smallmatrix}
  (-1)^{\epsilon_1} && 0 \\&&\\
  r && (-1)^{\epsilon_2}\end{smallmatrix}\Big)\,\Big|\,\,
  \epsilon_{1,2}\in\{0,1\}, r\in
  \IZ\Big\}\,,
 \ee
the following coset decomposition holds:
 \be\label{Dec12}
  GL_2(\IZ)=\bigsqcup_{(m,n)\in \CP}\,\,\g_{(m,n)}\,B(\IZ)\,,
 \ee
where $\CP$ is defined in \eqref{Range}. The matrices $\gamma_{(m,n)}$
are given by
 \be\label{gmn}
  \g_{(m,n)}=\Big(\begin{smallmatrix}
                    k && m\\&&\\l && n\end{smallmatrix}\Big),
                  \ee
where in \eqref{gmn} we choose  $n>0$ for $n\neq 0$ and $m>0$ for
$n=0$. The entries $k,l$ satisfy the equation
 \be
  \det\g_{(m,n)}=kn-lm=1\,,
 \ee
and are uniquely defined by the additional conditions
 \be
  0\leq l<n,\quad \,\,{\rm if}\,\,\,n\neq 0\,.
 \ee
For $n=0$ (and hence $m=1$) we take $k=0$ and $l=1$.
\end{lem}

\proof Given $\alpha=\Big(\begin{smallmatrix}
 k && m \\&&\\l && n\end{smallmatrix}\Big)\in GL_2(\IZ)$,
the right multiplication  by an  element
 \be
  \beta=\Big(\begin{smallmatrix}(-1)^{\e_1}&&0\\&&\\r&&(-1)^{\e_2}\end{smallmatrix}\Big)\in
  B(\IZ)\,,
 \ee
is given by
 \be
  \alpha\,\beta\,
  =\,\Big(\begin{smallmatrix}k&&m\\&&\\l &&n\end{smallmatrix}\Big)
  \Big(\begin{smallmatrix}(-1)^{\e_1}&&0\\&&\\r &&(-1)^{\e_2}\end{smallmatrix}\Big)
  =\Big(\begin{smallmatrix}
    mr+(-1)^{\e_1}k &&(-1)^{\e_2}m\\&&\\
  nr +(-1)^{\e_1}l&& (-1)^{\e_2}n\end{smallmatrix}\Big)\,.
 \ee
The determinant constraint
 \be
    |\det\alpha|\,=\,|kn-lm|=1\,,
 \ee
implies that
 \be
  {\rm gcd}(m;n)\,=1\,,
 \ee
and therefore we have a projection  on $\CP$
 \be\label{Col2}
  GL_2(\IZ) \longrightarrow\CP,\quad
  \Big(\begin{smallmatrix}k&&m\\&&\\l&&n\end{smallmatrix}\Big)\longmapsto(m,n)\in
  \CP\,.
 \ee
with the free action of $B(\IZ)$ on the fibers.
 We should demonstrate that the action on the fibers is transitive and
 pick a representative for each orbit. Thus for
 $n\neq 0$ we choose $n>0$ and $0\leq l<n$. Then, by the
Bezout identity, $k$ is uniquely determined from the determinant condition $\det
\alpha=1$. For $n=0$ we choose $m=1$ and thus one might
pick $k=0$ and $l=1$.  This proves the validity of the decomposition
\eqref{Dec12}.  $\Box$

\begin{cor}\label{Dec45}
Vector $\phi^{\vee}_{GL_2(\IZ)}$ in \eqref{SphVect2} allows the following
representation:
 \be\label{SphVect3}
  \phi^{\vee}_{GL_2(\IZ)}(x)\,
  =\!\sum_{\alpha\in GL_2(\IZ)/B(\IZ)}\!\!
  [\pi^{\vee}_{\g}(\alpha)\,\delta](x)\,.
 \ee
\end{cor}

\proof Applying \eqref{TrRules}, for
$\alpha=\Big(\begin{smallmatrix}k&&m\\&&\\l&&n\end{smallmatrix}\Big)\in
GL_2(\IZ)$,
 \be\label{TrRules1}
  [\pi^{\vee}_{\g}(\alpha)\,\delta](x)\,
   =\,|k+xl|^{\imath(\bar{\g}_1-\bar{\g}_2)-1}\,\delta\Big(\frac{m+xn}{k+xl}\Big)\\
  =\,|k+xl|^{\imath(\bar{\g}_1-\bar{\g}_2)}\,\delta(m+xn)\,
  =\,|n|^{-\imath(\bar{\g}_1-\bar{\g}_2)}\,\delta(m+xn)\,,
 \ee
due to substitution $x=-m/n$ into $k+xl$ and the fact that
$|\det\alpha|=1$. Therefore, by Lemma \ref{BZG} the summation over
$\CP$ is equivalent to the summation over $GL_2(\IZ)/B(\IZ)$. $\Box$

Now the Eisenstein automorphic function
associated with $(\pi_\gamma,\CV_\gamma)$ is defined as the following matrix  element
 \be \label{Autform}
  \Phi_{\gamma}(g)\,
  =\,\<\phi^{\vee}_{GL_2(\IZ)},\pi_\gamma(g) \phi_{O(2)}\>\,,
 \ee
where the pairing is defined in \eqref{PAIR}. 
Note that $\Phi_{\gamma}(g)$ obviously  defines a function on the
double coset $\CM_2$. Explicit realization  of the  principal series
representation allows to obtain an explicit expression for the
$GL_2$-automorphic function  \eqref{Autform}.

\begin{prop}
For the $GL_2$-Eisenstein function given by the matrix element \eqref{Autform}, 
the following series representation holds,
  \be\label{Eisen}
   \Phi_\gamma(g)\,
   =\,\<\phi^{\vee}_{GL_2(\IZ)}\,, \pi_\gamma(g)\, \phi_{O_2}\>\\
   =\,|\det g|^{\imath \gamma_2+\frac{1}{2}}\!
   \sum_{(n,m)\in\CP}\!
  \bigl|(na+mc)^2+(nb+md)^2\bigr|^{\frac{\imath(\gamma_1-\gamma_2)-1}{2}}\,, \quad
  g=\Big(\begin{smallmatrix}a&&b\\&&\\c&&d\end{smallmatrix}\Big)\,.
 \ee
provided ${\rm Im}(\g_1-\g_2)>1$ for convergence.
\end{prop}

\proof Substituting the expressions \eqref{SphVect1},
\eqref{SphVect2} for the vectors $\phi_{O_2}\,,\phi^{\vee}_{GL_2(\IZ)}$ in
\eqref{Autform} and considering \eqref{TrRules} we derive
 \be
  \Phi_{\gamma}(g)\,
  =\,\<\phi^{\vee}_{GL_2(\IZ)},\,\pi_\gamma(g)\,\phi_{O_2}\>\,
  =\!\int\limits_{\IR}\!dx\,\,
  \ov{\phi^{\vee}_{GL_2(\IZ)}(x)}\,\,[\pi_\gamma(g)\,\phi_{O_2}](x)\\
  =\,|\det(g)|^{\i\g_2+\frac{1}{2}}\!\sum_{(n,m)\in\CP}
  |n|^{\i(\g_1-\g_2)}\!
  \int\limits_{\IR}\!dx\,\,\delta(m+xn)\\
  \times\,\,\bigl|(a+xc)^2+(b+xd)^2\bigr|^{\frac{\i(\g_1-\g_2)-1}{2}}\\
  =\,|\det(g)|^{\i\g_2+\frac{1}{2}}\!\sum_{(n,m)\in\CP}
  |n|^{\i(\g_1-\g_2)-1}\,
  \Big|\Big(a-\frac{mc}{n}\Big)^2+\Big(b-\frac{md}{n}\Big)^2
  \Big|^{\frac{\i(\g_1-\g_2)-1}{2}}\,,
 \ee
which gives \eqref{Eisen} after changing the summation variable $m\mapsto-m$. $\Box$

The Eisenstein function \eqref{Autform} is right $O_2$-invariant and
hence is defined as a function on the coset space $GL_2(\IR)/O_2$.
We choose the following set of representatives of the right
$O_2$-cosets
 \be\label{Param123}
  g(\tau,t)\,=\,t^{\frac{1}{2}}\tau_2^{-\frac{1}{2}}\,
  \begin{pmatrix}
  \tau_2 & \tau_1\\
  0 &1\end{pmatrix}\in GL_2(\IR)/O_2\,, \quad \tau_1\in \IR, \quad
  t,\tau_2,\,t^{\frac{1}{2}},\tau_2^{\frac{1}{2}} \in \IR_+\,.
 \ee
Then evaluation  of the Eisenstein function \eqref{Eisen} on
elements \eqref{Param123} reads
 \be\label{EisenRed}
  \Phi_\gamma(\tau,t)\,
  =\,t^{\frac{\imath(\gamma_1+\gamma_2)}{2}}\!\!
  \sum_{(n,m)\in\CP}\,
  \frac{({\rm Im}(\tau))^{\frac{\i(\g_2-\g_1)+1}{2}}}
  {|m+n\tau|^{\i(\g_2-\g_1)+1}}\,,\qquad \tau=\tau_1+\imath
  \tau_2\,,
 \ee
provided by the direct calculation, for $g=g(\tau,t)$:
 \be
  |\det g|^{\imath\gamma_2+\frac{1}{2}}\,
  \bigl|(na+mc)^2+(nb+md)^2\bigr|^{\frac{\imath(\gamma_1-\gamma_2)-1}{2}}\\
  =\,t^{\imath\gamma_2+\frac{1}{2}}
  \Big|n^2t\tau_2
  +\frac{(nt^{\frac{1}{2}}\tau_1+mt^{\frac{1}{2}})^2}{\tau_2}
  \Big|^{\frac{\imath(\gamma_1-\gamma_2)-1}{2}}\\
  =\,t^{\frac{\imath(\g_1+\g_2)}{2}}\,\tau_2^{-\frac{\imath(\gamma_1-\gamma_2)-1}{2}}\,
  \bigl|m^2\tau_2^2+(m+n\tau_1)^2\bigr|^{\frac{\imath(\gamma_1-\gamma_2)-1}{2}}\,.
 \ee
Also the Eisenstein series may be written in a form that makes its
$GL_2(\IZ)$-invariance obvious.

\begin{prop} The $GL_2$-Eisenstein function \eqref{EisenRed} allows the
  following presentation
 \be \label{EisenInv}
  \Phi_\gamma(\tau,t)\,
  =\,t^{\frac{\imath(\gamma_1+\gamma_2)}{2}}\!\!
  \sum_{\alpha\in GL_2(\IZ)/B(\IZ)}
  ({\rm Im}(\alpha\cdot\tau))^{\frac{\i(\g_2-\g_1)+1}{2}}\,,
 \ee
where  the action of $GL_2(\IZ)$ is given by \eqref{RLaction}:
 \be
  \alpha\cdot\tau\,
  =\,\frac{k+\tau l}{m+\tau n},\qquad
  \alpha=\Big(\begin{smallmatrix}
  m&&k\\&&\\n&&l
  \end{smallmatrix}\Big)\in GL_2(\IZ)\,.
 \ee
\end{prop}

\proof The expression \eqref{EisenInv} follows from
\eqref{EisenRed}, for $\alpha=\Big(\begin{smallmatrix}
  m&&k\\&&\\n&&l
  \end{smallmatrix}\Big)$,
 \be
  {\rm Im}(\alpha\cdot\tau)\,
  =\,{\rm Im}\Big(\frac{k+\tau l}{m+\tau n}\Big)\,
  =\,\frac{1}{2\i}\Big(\frac{k+\tau l}{m+\tau n}
  -\frac{k+\bar{\tau}l}{m+\bar{\tau}n}\Big)\,
  =\,\frac{\det\a\,{\rm Im}(\tau)}{|m+\tau n|^2}\,.
 \ee
Alternatively, according to \eqref{SphVect3} we have
 \be
  \Phi_\gamma(\tau,t)\,
  =\!\!\sum_{\alpha\in GL_2(\IZ)/B(\IZ)}
  \bigl\<\pi^{\vee}_\gamma(\alpha)\,\delta\,,\,
  \pi_\gamma(g(\tau,t))\,\phi_{O_2}\bigr\>\,.
 \ee
Then considering the element \eqref{Param123} acting on the
spherical vector \eqref{SphVect1} via \eqref{TrRules} gives
 \be
  \<\delta\,,\pi_\gamma(g(\tau,t))\,\phi_{O_2}\>\,
  =\int\limits_{\IR}\!dx\,\,\delta(x)\,
  [\pi_\gamma(g(\tau,t))\,\phi_{O_2}](x)\\
  =\,\int\limits_{\IR}\!dx\,\delta(x)\,t^{\imath\g_2+\frac{1}{2}}
  |t\tau_2|^{\frac{\imath(\g_1-\g_2)-1}{2}}
  \left(1+\Big(\frac{\tau_1\sqrt{\frac{t}{\tau_2}}+x\sqrt{\frac{t}{\tau_2}}}
  {\sqrt{t\tau_2}}\Big)^2\right)^{\frac{\i(\g_1-\g_2)-1}{2}}\\
  =\,t^{\frac{\imath(\gamma_1+\gamma_2)}{2}}
  \,|\tau_1^2+\tau_2^2|^{\frac{\i(\g_2-\g_1)+1}{2}}\,
  =\,t^{\frac{\imath(\gamma_1+\gamma_2)}{2}}
  \,|\tau|^{\i(\g_2-\g_1)+1}\,,
 \ee
which entails the presentation \eqref{EisenInv}.  $\Box$

\section{ Global $GL_2$ Hecke-Baxter operator}

Global $GL_2$ Hecke-Baxter operator is represented by an element of
the associative algebra
$\CH(GL_2(\IQ),GL_2(\IZ))\otimes\CH(GL_2(\IR),O_2)$ acting via
convolutions (from the left and from the right) in the space of
functions on $\CM_2=GL_2(\IZ)\backslash GL_2(\IR)/O_2$. To construct
such an operator we start with
 describing elements of the Hecke algebras
$\CH(GL_2(\IQ),GL_2(\IZ))$  and $\CH(GL_2(\IR),O_2)$ acting in the
space of the $GL_2$-automorphic forms.

Taking into account the interpretation of  $\CM_2$ as a space of
equivalence classes of lattices $L$ in $\IR^2$ we introduce the
following  averaging operators $T_n,\,n\in\IZ_+$ as analogs of
\eqref{Tn} acting via
 \be\label{HA1}
  (T_n\cdot f)(L)\,
  =\!\sum_{[L:L']=n}f(L')\,,
 \ee
where the sum goes over sub-lattices $L'\subset L$ of index $n$.
The double coset description of $\CM_2$ allows to rewrite
the action of operators \eqref{HA1} as follows
 \be\label{HA2}
  (T_n\cdot f)(g)\,
  =\!\!\sum_{\gamma\in {GL_2(\IZ)\backslash\Mat}^{(n)}_2(\IZ)}\!\!
  f(\gamma\, g)\,,
 \ee
where
 \be\label{Matn}
  \Mat^{(n)}_2(\IZ)\,
  =\,\bigl\{\gamma\in {\rm Mat}_2(\IZ)\,\bigr|\quad|\det \gamma|=n\bigr\}\,.
 \ee
This action may be written more explicitly using a specific choice
of coset representatives of the quotient space
$GL_2(\IZ)\backslash\Mat^{(n)}_2(\IZ)$.

\begin{lem}
For each $n\in\IZ_+$, the space $\Mat^{(n)}_2(\IZ)$ defined in
  \eqref{Matn} allows
  the following coset decomposition:
 \be\label{Coset1}
  \Mat^{(n)}_2(\IZ)\,
  =\bigsqcup_{i=1}^{\s(n)}\,\,GL_2(\IZ)\,\alpha_i\,,
 \ee
where $\s(n)=\sum_{d|n}d$, and
 \be\label{Cosetrep}
  \alpha_i=\begin{pmatrix} a_i& b_i\\0 & d_i\end{pmatrix}, \qquad
  a_i,b_i,d_i>0,\quad a_id_i=n, \quad 0\leq b_i<d_i\,.
 \ee
In particular, for each $x\in\Mat_2^{(n)}(\IZ)$ there exists a
unique $\gamma\in GL_2(\IZ)$ such that
 \be\label{Dec2}
  x\,=\,\gamma\,\alpha_i\,,
 \ee
for some $\alpha_i$ of the form \eqref{Cosetrep}.
\end{lem}

\proof  To find proper set of coset representatives
we first check  that
the  left multiplication by a matrix $\Lambda \in
GL_2(\IZ)$ allows us to put the lower left element of the matrix to zero.
We have
 \be
  \begin{pmatrix}
  \lambda_1 & \lambda_2\\
  \lambda_3 & \lambda_4\end{pmatrix}
  \begin{pmatrix}
  a & b\\
  c& d \end{pmatrix}\,
  =\,
  \begin{pmatrix}
  * & *\\
  a \lambda_3+c\lambda_4  & *\end{pmatrix}\,.
 \ee
Thus we should find $(\lambda_3,\lambda_4)$ such that
 \be
  a\lambda_3+c\lambda_4=0\,.
 \ee
We fullfill this condition by taking $\lambda_4=a/{\rm gcd}(a;c)$,
$\lambda_3=-c/{\rm gcd}(a;c)$, so that such $\la_3,\la_4$ are
relatively prime and hence represent a row of an invertible matrix
$\Lambda$. Thus we have
 \be
  \begin{pmatrix}
  \lambda_1 & \lambda_2\\
  \lambda_3&\lambda_4\end{pmatrix}
  \begin{pmatrix}
  a & b\\
  0& d \end{pmatrix}\,
  =\,
  \begin{pmatrix}
  a\lambda_1 & b\lambda_1+d\lambda_2\\
  a\lambda_3  & b\lambda_3+d\lambda_4\end{pmatrix}\,.
 \ee
To retain the constraint $c=0$ we should have $a\lambda_3=0$. As $a\neq
0$ (since  $ad\neq 0$), hence we have $\lambda_3=0$  thus arriving  at
 \be
  \begin{pmatrix} \lambda_1 & \lambda_2\\ 0&
                                                  \lambda_4\end{pmatrix}
  \begin{pmatrix} a & b\\ 0& d \end{pmatrix}=
  \begin{pmatrix} a\lambda_1 & b\lambda_1+d\lambda_2
  \\  0  & d\lambda_4\end{pmatrix}\,,
 \ee
where
 \be
  |\lambda_1|=|\lambda_4|=1\,.
 \ee
Therefore we might put $a,d\in \IZ_+$ so that what  remains are the
following  transformations
 \be
  \begin{pmatrix}
  1 & \lambda_2\\
  0&1 \end{pmatrix}
  \begin{pmatrix}
  a & b\\
  0& d \end{pmatrix}\,
  =\,
  \begin{pmatrix}
  a & b+d\lambda_2\\
  0  & d\end{pmatrix}\,.
 \ee
This allows us to make  $b$ to satisfy the condition $0\leq b<d$.
This proves \eqref{Coset1}.  Direct check using
 \be
  \begin{pmatrix}
  \lambda_1 & \lambda_2\\
  \lambda_3&\lambda_4\end{pmatrix}
  \begin{pmatrix}
  a & b\\
  0& d \end{pmatrix}\,
  =\,
  \begin{pmatrix}
  a\lambda_1 & b\lambda_1+d\lambda_2\\
  a\lambda_3  & b\lambda_3+d\lambda_4\end{pmatrix}\,,
 \ee
shows that the stabilizer  of such elements is trivial. $\Box$

As a consequence of the previous Lemma we obtain the following
presentation for the operators $T_n$ in \eqref{HA2} acting via
 \be
  (T_n\cdot f)(g)\,
  =\sum_{a,d>0\atop ad=n}\sum_{b=0}^{d-1}
  f\Big(\,\Big(\begin{smallmatrix} a&&b\\&&\\0&&d\end{smallmatrix} \Big)\, g\Big)\,.
 \ee
Let us combine the operators $T_n$, $n\in \IZ_+$ into the  generating series
 \be
  Q^{GL_2(\IZ)}_s\,
  =\,\sum_{n=1}^{\infty}\, \frac{1}{n^{s+\frac{1}{2}}}\, T_n\,,
 \ee
acting on a function in the space of lattices $L\subset\IR^2$ in the
following way
 \be
  (Q^{GL_2(\IZ)}_s\cdot f)(L)\,
  =\sum_{L'\subset L}^{\infty}\, \frac{1}{[L:L']^{s+\frac{1}{2}}}\,f(L')\,.
  \ee
Equivalently,
 \be\label{HBZ0}
  \bigl(Q^{GL_2(\IZ)}_s\cdot f\bigr)(g)\,
  =\!\!\sum_{\a\in GL_2(\IZ)\backslash\Mat^*_2(\IZ)}\,\,
  \frac{1}{|\det\a|^{s+\frac{1}{2}}} \,f(\alpha\, g)\,,
 \ee
where $\a$ runs through the set of $GL_2(\IZ)$-coset representatives
\eqref{Cosetrep} of the space
 \be
  \Mat^*_2(\IZ)\,=\,\Mat_2(\IZ)\cap GL_2(\IQ)\,
  =\bigsqcup_{n\in \IZ_+} \Mat^{(n)}_2(\IZ)\,.
 \ee
Introduce the following  analog of the Riemann zeta-function for
$GL_2$:
 \be\label{Riem}
  \zeta^{GL_2}(s)\,
  =\!\!\sum_{\alpha\in  GL_2(\IZ)\backslash {\rm Mat}^*_2(\IZ)}\,
  \frac{1}{|\det  \alpha|^{s+\frac{1}{2}}}\,,\quad{\rm Re}(s)>\frac{3}{2}\,.
 \ee

\begin{lem} For the function \eqref{Riem}, the following identity
holds:
 \be\label{2Riem}
  \zeta^{GL_2}(s)\,
  =\!\!\sum_{\alpha\in GL_2(\IZ)\backslash {\rm Mat}^*_2(\IZ)}\,
  \frac{1}{|\det  \alpha|^{s+\frac{1}{2}}}\,
  =\,\zeta\Big(s+\frac{1}{2}\Big)\,
  \zeta\Big(s-\frac{1}{2}\Big)\,,
 \ee
where $\zeta(s)$ is the Riemann zeta-function given by
 \be
  \zeta(s)\,
  =\,\sum_{n=1}^{\infty} \,\frac{1}{n^s}\,, \qquad {\rm Re}(s)>1\,.
 \ee
\end{lem}

\proof Using the set \eqref{Cosetrep} of  coset representatives
we find out
 \be
  \zeta^{GL_2}(s)\,
  =\!\!\sum_{\alpha \in GL_2(\IZ)\backslash {\rm Mat}^*_2(\IZ)}\,
  \frac{1}{|\det  \alpha |^{s+\frac{1}{2}}}\,
  =\sum_{a,d>0 \atop 0\leq b<d}
  \frac{1}{a^{s+\frac{1}{2}}}\,
  \frac{1}{d^{s+\frac{1}{2}}}\,.
 \ee
Summation over $b$ results in
 \be
  \sum_{\alpha \in GL_2(\IZ)\backslash {\rm Mat}^*_2(\IZ)}\,
  \frac{1}{|\det  \alpha|^{s+\frac{1}{2}}}\,
  =\sum_{a,d>0}\frac{1}{a^{s+\frac{1}{2}}}\,
  \frac{1}{d^{s-\frac{1}{2}}}=
  \zeta\Big(s+\frac{1}{2}\Big)\,\zeta\Big(s-\frac{1}{2}\Big)\,,
 \ee
which completes the proof. $\Box$

Below we encounter a generalization $\zeta^{GL_2}(s|\gamma)$ of the
function \eqref{2Riem} associated with a spherical principal series
representation  $(\pi_\gamma,\CV_\gamma)$ of $GL_2(\IR)$.

\begin{prop}\label{HBZ}
The action of the operator $Q_s^{GL_2(\IZ)}$ on the
$GL_2$-Eisenstein functions \eqref{Eisen} associated with a
spherical principal series representation $(\pi_\gamma,\CV_\gamma)$
is given by
 \be \label{Act456}
  \bigl(Q_s^{GL_2(\IZ)}\cdot \Phi_\gamma\bigr)(g)\,
  =\,\zeta^{GL_2}(s|\gamma)\,\,\Phi_\gamma(g)\,,
 \ee
where
 \be\label{CompGLzeta}
 \zeta^{GL_2}(s|\gamma)=\zeta(s-\imath\g_1)\,
 \zeta(s-\imath\g_2)\,.
\ee
\end{prop}

\proof
 Application of the Hecke-Baxter operator \eqref{HBZ0} to the matrix element \eqref{Autform} reads
 \be
  \bigl(Q_s^{GL_2(\IZ)}\cdot \Phi_\gamma\bigr)(g)
  =\sum_{\alpha \in GL_2(\IZ)
  \backslash {\rm Mat}^*_2(\IZ)}\,
  \frac{1} {|\det \alpha |^{s+\frac{1}{2}}}\,\,
  \Phi_\gamma(\alpha \cdot g)\,,
 \ee
where  
 \be
  \Phi_\gamma(\alpha \cdot g)\,
  =\,\bigl\<\phi^{\vee}_{GL_2(\IZ)}\,,
  \pi_{\gamma}(\alpha\,g)\,\phi_{O_2}\bigr\>\,
  =\,\bigl\<\pi^{\vee}_{\gamma}(\alpha^{-1})\,\phi^{\vee}_{GL_2(\IZ)}\,,
  \pi_{\gamma}(g)\,\phi_{O_2}\bigr\>\,.
 \ee
Below we calculate the $Q_s^{GL_2(\IZ)}$-action on the
$GL_2(\IZ)$-invariant vector \eqref{SphVect2}:
 \be
  \bigl(Q_s^{GL_2(\IZ)}\cdot\phi^{\vee}_{GL_2(\IZ)})(x)\,
  =\!\!\sum_{\alpha \in GL_2(\IZ)\backslash {\rm Mat}^*_2(\IZ)}\,
  \frac{1} {|\det \alpha|^{s+\frac{1}{2}}}\,\,
  \bigl[\pi^{\vee}_{\gamma}(\alpha^{-1})\,\phi^{\vee}_{GL_2(\IZ)}\bigr](x)\,.
 \ee
Using the presentation \eqref{SphVect3} for the $GL_2(\IZ)$-invariant
vector,
 \be\label{SphVect34}
  \phi^{\vee}_{GL_2(\IZ)}(x)\,
  =\!\sum_{\beta\in GL_2(\IZ)/B(\IZ)}\![\pi^{\vee}_{\g}(\beta)\,\delta](x)\,,
 \ee
we deduce the following:
 \be\label{Qact}
  \bigl(Q_s^{GL_2(\IZ)}\cdot \phi^{\vee}_{GL_2(\IZ)})(x)\\
  =\!\sum_{\alpha \in GL_2(\IZ)\backslash\Mat^*_2(\IZ)}\,
  \frac{1} {|\det \alpha |^{s+\frac{1}{2}}}\,\,
  \sum_{\beta\in GL_2(\IZ)/B(\IZ)}\,\,
  [\pi^{\vee}_{\gamma}(\alpha^{-1}\beta)\,\delta](x)\\
  =\!\sum_{\alpha\in B(\IZ)\backslash\Mat^*_2(\IZ)}\,
  \frac{1} {|\det\alpha|^{s+\frac{1}{2}}}\,\,
  [\pi^{\vee}_{\gamma}(\alpha^{-1})\,\delta](x)\,.
 \ee
Given $\alpha=\Big(\begin{smallmatrix} k&&
m\\&&\\l&& n\end{smallmatrix}\Big)\in{\rm Mat}^*_2(\IZ)$, consider
the adjugate matrix $\beta\in\Mat_2^*(\IZ)$:
 \be
  \beta=\Big(\begin{smallmatrix} 
  n && -m\\&&\\-l && k\end{smallmatrix}\Big)\in {\rm Mat}^*_2(\IZ)\,,\quad 
  \alpha^{-1}=(\det\beta)^{-1}\beta, \qquad \det \alpha=\det \beta\,. 
 \ee
 Then substituting into \eqref{Qact} gives 
 \be\label{Sum12}
  \bigl(Q_s^{GL_2(\IZ)}\cdot\phi^{\vee}_{GL_2(\IZ)})(x)
   =\!\sum_{\beta \in  {\rm Mat}^*_2(\IZ)/B(\IZ)}\,
 \,\, \frac{1} {|\det \beta|^{s+\frac{1}{2}}}\,\,
  [\pi^{\vee}_{\gamma}\bigl((\det\beta)^{-1}\beta\bigr)\,\delta](x)\,.
  \ee
In the following we
  change the notations and use the following (equivalent) parameterization for 
$\beta=\Big(\begin{smallmatrix} k&&
m\\&&\\l&& n\end{smallmatrix}\Big)\in{\rm Mat}^*_2(\IZ)$
in \eqref{Sum12}.   
Recall the explicit form for the representation action
\eqref{TrRules1}, for thus defined parmaterization of the
elements of ${\rm Mat}^*_2(\IZ)$: 
 \be\label{DeltaACT}
  [\pi^{\vee}_\gamma(\beta)\,\delta](x)\,
  =\,|\det \beta|^{\imath \bar{\gamma_2}+\frac{1}{2}}
  |k+xl|^{\imath (\bar{\gamma_1}-\bar{\gamma_2})-1}\,
  \delta\left(\frac{m+ xn}{k+xl}\right)\\
  =\,|\det \beta|^{\imath \bar{\gamma_1}+\frac{1}{2}}\,
  |n|^{\imath (\bar{\gamma_2}-\bar{\gamma_1})}\,\delta(m+ xn)\,.
 \ee
Then we introduce
 \be
  p={\rm gcd}(m,n)>0,\qquad q={\rm gcd}(k,l)>0\,,
 \ee
and rewrite the element $\alpha\in{\rm Mat}^*_2(\IZ)$ as follows:
 \be
  \alpha\,
  =\,\Big(\begin{smallmatrix}
  qk'&&pm'\\&&\\ql'&&pn'\end{smallmatrix}\Big)\,
  =\,\Big(\begin{smallmatrix} 
  k'&& m'\\&&\\l'&& n'\end{smallmatrix}\Big)
  \Big(\begin{smallmatrix}
  q&&0\\&&\\0&&p\end{smallmatrix}\Big)\,, \qquad
  {\rm gcd}(m',n')={\rm gcd}(k',l')=1\,,
 \ee
providing
 \be\label{alpha}
  \beta\,
  =\,\Big(\begin{smallmatrix}
  k&&m\\&&\\l&& n\end{smallmatrix}\Big)\,
  =\,g\,\Big(\begin{smallmatrix}
  q&&0\\&&\\0&&p\end{smallmatrix}\Big)\,,\quad
  g\in GL_2(\IZ),\,\quad
  \det\beta\,=\,pq\,.
  \ee
Hence we have 
 \be\label{beta}
  (\det \beta)^{-1}\,\beta\,
  =\,g\,\Big(\begin{smallmatrix}
  p^{-1}&&0\\&&\\0&&q^{-1}\end{smallmatrix}\Big)\,.
  \ee
  
Now for a fixed column $(m,n)$ in \eqref{alpha}, there is a unique
pair of $(k,l)$ (up to the $B(\IZ)$-action), such that
$\det \alpha=pq$. However there is subtlety here: taking a
quotient over $B(\IZ)$ we may use only the subgroup $B(\IZ)^p\subset
B(\IZ)$ (to transform $k,l$),
 \be
  B(\IZ)^p\,
  =\,\Big\{\Big(\begin{smallmatrix}
  (-1)^{\e_1}&&0\\&&\\
  pr&&(-1)^{\e_2}
  \end{smallmatrix}\Big):\,\,\e_{1,2}\in\{0,1\},\,r\in\IZ\Big\}\,,\quad
  |B(\IZ)/B(\IZ)^p|=p\,,
 \ee
hence  we acquire an additional factor $p$.
By \eqref{DeltaACT}, for the diagonal matrix in \eqref{beta} we have   
 \be
  [\pi^{\vee}_\gamma({\rm diag}(p^{-1},q^{-1})\,\delta](x)\,
  =(pq)^{-\imath\bar{\gamma_1}-\frac{1}{2}}\, q^{\imath
    (\bar{\gamma_1}-\bar{\gamma_2})+1}\,\,\,\delta(x)\\
    =\,p^{-\imath\bar{\gamma_1}-\frac{1}{2}}\, q^{-\imath
  \bar{\gamma_2}+\frac{1}{2}}\,\,\,\delta(x)\,. 
 \ee
Now the action of the operator $Q_s^{GL_2(\IZ)}$ takes the following
form 
 \be
  \bigl(Q_s^{GL_2(\IZ)}\cdot \phi^{\vee}_{GL_2(\IZ)}\bigr)(x)
  =|B(\IZ)/B(\IZ)^p|\,
  \Big(\sum_{p,q>0}
  \frac{p^{-\imath\bar{\gamma_1}-\frac{1}{2}}\, q^{-\imath
  \bar{\gamma_2}+\frac{1}{2}}}{(pq)^{s+\frac{1}{2}}}\Big)\,
  \phi^{\vee}_{GL_2(\IZ)}(x)\,.
 \ee
Thus we finally arrive at
 \be
  \bigl(Q_s^{GL_2(\IZ)}\cdot
  \phi^{\vee}_{GL_2(\IZ)}\bigr)(x)\,
  =\,\zeta(s+\imath \bar{\gamma}_1)\,\zeta(s+\imath
  \bar{\gamma}_2)\,\,\phi^{\vee}_{GL_2(\IZ)}(x)\,. 
  \ee
  Taking into account that in the matrix elements the left vector is
  taken complex conjugated  we arrive at \eqref{Act456}. $\Box$


The Archimedean counterpart of the one-parameter family
$Q^{GL_2(\IZ)}_s$   of elements in \\$\CH(GL_2(\IQ),GL_2(\IZ))$
given by \eqref{HBZ0} is the
 $GL_2(\IR)$ Hecke-Baxter operator \cite{GLO08}:
 \be
  \bigl(Q^{GL_2(\IR)}_s\cdot f\bigr(\tilde{g})\,
  =\!\int\limits_{GL_2(\IR)}\!\! d\mu^G(g)\,\, |\det g|^{s+\frac{1}{2}} f(g^{-1}\tilde{g})\,,\\
  d\mu^G(g)=e^{-\pi \Tr(g^{\top}g)}\,d\mu(g)\,,
 \ee
acting via convolution with the following $O_2$-biinvariant function
on $GL_2(\IR)$:
 \be\label{ArKernel2}
  Q^{GL_2(\IR)}_s(g)=|\det g|^{s+\frac{1}{2}}\,e^{-\pi \Tr(g^{\top}g)}\,.
 \ee
Its action on the matrix element \eqref{Autform} (actually on any
matrix element with the right $O_2$-invariant vector) was calculated in
\cite{GLO08} and is given by
 \be\label{Act123}
  \bigl(Q^{GL_2(\IR)}_s\cdot \Phi_\gamma\bigr)(g)\,
  =\,L^{GL_2(\IR)}(s|\gamma)\,
  \Phi_\gamma(g)
 \ee
where
 \be\label{GammaF}
  L^{GL_2(\IR)}(s|\gamma)=
  \prod_{j=1}^2 \pi^{-\frac{s-\imath \gamma_j}{2}}\,
  \Gamma\Big(\frac{s-\imath \gamma_j}{2}\Big)\,.
 \ee

Let us note that while the operator $Q_s^{GL_2(\IZ)}$ acts by
convolution on the functions on the double coset space
$\CM_2=GL_2(\IZ)\backslash GL_2(\IR)/O_2$ from the left, the
operator $Q_s^{O_2}$ acts by convolution on the functions on the
space $\CM_2$ from  the right. Still it is reasonable to consider
its combination acting by simultaneous left/right convolution. Let
us  define  the following integral operator acting on
$GL_2$-automorphic functions
 \be\label{GHBO1}
  \bigl(\wh{Q}^{GL_2}_s \bullet  \Phi\bigr)(\tilde{g})\,
  =\!\int\limits_{GL_2(\IR)}\!\!d\mu(g)\,\,\wh{Q}^{GL_2}_s(\tilde{g},g)\,\,\Phi(g^{-1})\,,
 \ee
where the kernel of the integral operator is given by
 \be\label{GHBO}
  \wh{Q}^{GL_2}_s(\tilde{g},g)=\,\,|\det(\tilde{g}g)|^{s+\frac{1}{2}}\!\!
  \sum_{\a\in GL_2(\IZ)\backslash {\rm Mat}^*_2(\IZ)}\!\!
  e^{-\pi\Tr(g^{\top}g
  \,\alpha\,\tilde{g}\tilde{g}^{\top}\,\alpha^{\top})}\,,
 \ee
where the summation goes over representatives \eqref{Cosetrep}. One
might rewrite this operator in the following form:
 \be\label{GHBO1x}
  \bigl(\wh{Q}^{GL_2}_s \bullet  \Phi\bigr)(\tilde{g})\,
  =\!\!\int\limits_{GL_2(\IR)/GL_2(\IZ)}\!\!\!\!d\mu(g)\,\,
  \tilde{Q}^{GL_2}_s(\tilde{g},g)\,\,\Phi(g^{-1})\,,
 \ee
with the following integral kernel
 \be\label{GHBOx}
  \tilde{Q}^{GL_2}_s(\tilde{g},g)\,
  =\,|\det(\tilde{g}g)|^{s+\frac{1}{2}}\!\!
  \sum_{\a\in {\rm Mat}^*_2(\IZ)}\!\!
  e^{-\pi\Tr(g^{\top}g\,\alpha\,\tilde{g}\tilde{g}^{\top}\,\alpha^{\top})}\,.
\ee

This operator is a global analog  (in the sense of arithmetic geometry
of $\overline{{\rm Spec}(\IZ)}$) of the local
 spherical Hecke-Baxter  operators considered above. Our previous
 considerations may be summarized in the following form.

\begin{te}\label{TH1}  The global Hecke-Baxter operator \eqref{GHBO1},\eqref{GHBO1x}
  acts on the Eisenstein functions represented by matrix elements
\eqref{Eisen} of the spherical principal series representation
$(\pi_\gamma,\CV_\gamma)$ of $GL_2(\IR)$ by multiplication on the
corresponding global completed zeta-functions
 \be\label{GlobalL}
  \xi^{GL_2}(s|\gamma)\,
  =\prod_{j=1}^2\pi^{-\frac{s-\imath \gamma_j}{2}}\,
  \Gamma\Big(\frac{s-\imath \gamma_j}{2}\Big)\,
  \zeta(s-\imath \gamma_j)\,.
 \ee
\end{te}

\proof The action of $\wh{Q}^{GL_2}_s$ on \eqref{Eisen} is given by
the following integral
 \be
  \bigl(\wh{Q}^{GL_2}_s \bullet  \Phi_\gamma\bigr)(\tilde{g})\,
  =\,|\det\tilde{g}|^{s+\frac{1}{2}}\\
  \times\!\!\!\! \sum_{\a\in GL_2(\IZ)\backslash {\rm Mat}^*_2(\IZ)}\,\,
  \int\limits_{GL_2(\IR)}\!\!
  d\mu(g)\,\,|\det g|^{s+\frac{1}{2}} \,
  e^{-\pi\Tr[g\,\alpha\,\tilde{g}(g\alpha\,\tilde{g})^{\top}]} \,\,
  \Phi_\gamma(g^{-1})\,.
 \ee
Changing the integration variable $g$ by $g':=g\alpha\tilde{g}$ we obtain
 \be
  g=g'\tilde{g}^{-1}\alpha^{-1}\,,\quad
  g^{-1}=\alpha\,\tilde{g}(g')^{-1}\,,
 \ee
which results in (keeping the same notation $g$ for the integration
variable),
 \be
  \bigl(\wh{Q}^{GL_2}_s \bullet  \Phi_\gamma\bigr)(\tilde{g})\,
  =\,|\det\tilde{g}|^{s+\frac{1}{2}}\\
  \times\!\!\sum_{\a\in GL_2(\IZ)\backslash\Mat^*_2(\IZ)}\,\,
  \int\limits_{GL_2(\IR)}\!\!d\mu(g)\,\,
  \frac{|\det g|^{s+\frac{1}{2}}}{|\det\alpha|^{s+\frac{1}{2}}}\,\,
  e^{-\pi\Tr(g^{\top}g)}\,\,\Phi_\gamma(\alpha\tilde{g}g^{-1})\,.
 \ee
The above expression is the composition of the operators
$Q_s^{GL_2(\IZ)}$ and $Q_s^{GL_2(\IR)}$. Therefore, the assertion
follows from \eqref{Act456} and \eqref{Act123}. $\Box$

There is a global analog of the representation \eqref{ZetaInt} that
may be described as follows. Let us define  the $GL_2$-analogs of
the (modified) theta-constants \eqref{Tconst},  \eqref{Tconst12},
\eqref{Tconst123}:
 \be\label{Theta201}
  \Theta(0|T)\,
  =\sum_{\alpha \in\Mat_2(\IZ)}
    e^{\imath \pi\Tr(T\alpha\alpha^\top)}\,,
 \ee
 \be\label{Theta2}
  \Theta^*(0|T)\,
  =\sum_{\alpha \in\Mat_2^*(\IZ)}
    e^{\imath \pi\Tr(T\alpha\alpha^\top)}\,,
 \ee
 \be\label{Theta23}
  \Theta^{**}(0|T)\,
  =\sum_{\alpha \in  GL_2(\IZ)\backslash\Mat_2^*(\IZ)}
    e^{\imath \pi\Tr(T \alpha\alpha^\top)}\,,
 \ee
where $T$ is a complex symmetric matrix with a positive-definite
 imaginary part. Note that the series \eqref{Theta201} satisfies the following duality relation:
 \be\label{ThetaId}
  \Theta(0|T)=(\det(-\imath T))^{1/2}\, \Theta(0|-T^{-1})\,.
 \ee
Then the  global zeta-function \eqref{GlobalL}
 for $\gamma=0$ allows the following integral expression
  \be\label{ALL}
   \xi^{GL_2}(s|0)\,
   =\,\pi^{-s}\,\Gamma\Big(\frac{s}{2}\Big)^2\,\zeta(s)^2\\
   =\!\int\limits_{GL_2(\IR)/GL_2(\IZ)}\!\!
   d\mu(g)\,\,|\det g|^s\,
   \Theta^*(0|\imath g^{\top}g)\,.
  \ee
Indeed substituting  \eqref{Theta2} into \eqref{ALL} and making the
change of integration variable $h\to \alpha^{-1}h$ leads to
factorization of summation and integration. Thus  \eqref{ALL}
reduces to the product of zeta-functions and Gamma-factor
\eqref{GammaF} with $\gamma=0$. Note that the essential part of the
integral kernel \eqref{GHBO} is expressed in terms of the following
generalization of the classical theta constant
 \be
  \hat{\Theta}(0|A,B)\,
  =\!\sum_{\a\in {\rm  Mat}_2(\IZ)}\!\!\!\!e^{-\pi\Tr(A\alpha
    B\alpha^{\top})}\,.
 \ee
where $A$ and $B$ are symmetric positive-definite
$(2\times2)$-matrices. This kind of theta series is instrumental for
the verification of the analog of the functional equation for the
global $GL_2$ Hecke-Baxter operator extending the relations
\eqref{FR1} for $GL_1$. The corresponding functional equation for
the global Hecke-Baxter operator are compatible with the functional
equations for the completed  $GL_2$ zeta-function
 \be\label{GlobalL123}
  \xi^{GL_2}(1-s|-\gamma)\,=\,\xi^{GL_2}(s|\gamma)\,.
 \ee
In the Appendix we provide a proof of the functional
relation \eqref{GlobalL123} for $\gamma=0$.

Let us stress that the analytic functional relation
\eqref{GlobalL123} with $\g=0$ follows trivially from the explicit
expression \eqref{GlobalL}. The relevance of the proof given in
Appendix is in   explicit calculations  of various matrix integral
contributions mimicking the Riemann proof for $GL_1$ zeta-function.
This proof essentially uses the  transformation properties  of the
(incomplete) theta series \eqref{Theta2} and \eqref{Theta23}
entering the description of the global Hecke-Baxter operator and
thus demonstrates compatibility of the global Hecke-Baxter operator
$\wh{Q}_s^{GL_2}$ with the functional equation \eqref{GlobalL123}.
To complete the picture let us note that similar functional equation
hold for the $GL_2$-Eisenstein functions $\Phi_{\g}(g)$.

As a final remark let us observe that in the case of $GL_2$
(similarly to the case of $GL_1$, see  Section 2) one might
construct an interpolation of the global and Archimedean
Hecke-Baxter operators by considering congruence semigroups
 \be
  \Gamma(N)\,
  =\,\Id\,+\,N\Mat_2(\IZ)\subset GL_2(\IZ)\,.
 \ee
This might be useful to compare this approach with the results by D.
Kazhdan  \cite{Kaj}.


\section{Appendix:  Functional equation for $\xi^{GL_2}(s)$}

The aim of the Appendix is to prove the functional relation
 \be\label{GlobalL1234}
  \xi^{GL_2}(1-s|0)\,=\,\xi^{GL_2}(s|0)\,,
 \ee
satisfied by the function $\xi^{GL_2}(s|0)$ defined in
\eqref{GlobalL} using the integral representation \eqref{ALL}:
 \be\label{ALL123}
  \xi^{GL_2}(s|0)\,=\!\int\limits_{GL_2(\IR)/GL_2(\IZ)}\!\!
  d\mu(g)\,\,|\det g|^s\,
  \Theta^*(0|\imath g^{\top}g)\,.
 \ee
Here $d\mu(g)$ is the (restriction of the) Haar measure on
$GL_2(\IR)/GL_2(\IZ)$ and the theta-constant $\Theta^*(0|T)$ is
defined in \eqref{Theta201}. This proof generalizes the Riemann (second)
proof of the functional relation for the completed Riemann
zeta-function $\xi(s)$ to the case of $GL_2$.

We start with introducing the following subset of $(2\times
2)$-matrices
 \be
  \Mat_2^{(0)}(\IZ)\,
  =\,\{\gamma\in {\rm  Mat}_2(\IZ)\,|\quad\det \gamma=0\}\,.
 \ee

\begin{lem} The following decomposition holds
 \be
  \Mat_2^{(0)}(\IZ)\,
  =\,\Big(\bigsqcup_{\beta_I}\,GL_2(\IZ)\,\beta_{I}\Big)\,
  \cup\,
  \Big(\bigsqcup_{\beta_{II}}\,GL_2(\IZ)\,\beta_{II}\Big)\,,
 \ee
where the coset representatives are given by
 \be\label{Coset3}
  \beta_{I}=\begin{pmatrix} a & b\\0& 0\end{pmatrix}, \qquad a\in \IZ_+\,\quad b\in \IZ\,,
 \ee
and
 \be\label{Coset0}
  \beta_{II}=\begin{pmatrix} 0 & 0\\ 0&1\end{pmatrix}\,,
 \ee
with the corresponding  stabilizer subgroups of $GL_2(\IZ)$ being
 \be\label{Stab}
 {\rm St}_{I}=\tilde{B}_+(\IZ)\,
 =\,\Big\{ \Big(\begin{smallmatrix}
 \pm 1&& m\\&&\\0 &&1\end{smallmatrix}\Big)\Big|\,m\in \IZ\Big\}\,, \\
 {\rm St}_{II}=\tilde{B}_-(\IZ)\,
 =\,\Big\{ \Big(\begin{smallmatrix}
 \pm 1&& 0\\&&\\m &&1\end{smallmatrix}\Big)\Big|\,m\in \IZ\Big\}\,.
 \ee
\end{lem}

\proof  To find a proper set of coset representatives we first check
that left multiplication by a matrix $\Lambda \in GL_2(\IZ)$ allows
to put the lower left element of the matrix to zero. Indeed, we have
 \be
  \Lambda\,\g\,
  =\,\begin{pmatrix} \lambda_1 & \lambda_2\\\lambda_3&
                                                  \lambda_4\end{pmatrix}
  \begin{pmatrix} a & b\\ c& d \end{pmatrix}=
  \begin{pmatrix} * & *\\ a \lambda_3+c\lambda_4  &
  *\end{pmatrix}\,,
 \ee
we we need to find $(\lambda_3,\lambda_4)$ such that
 \be
  a \lambda_3+c\lambda_4=0\,.
 \ee
It is easy to fulfill this condition by
 taking   $\lambda_4=a/{\rm gcd}(a;c)$ and
$\lambda_3=-c/{\rm gcd}(a;c)$.  The  condition $c=0$ implies
$a\lambda_3=0$ allowing two solutions $a=0$ and $a\neq 0$,
$\lambda_3=0$.

Let us first consider the case of $a=0$, then we have
 \be
  \Lambda\,\g\,
  =\,\begin{pmatrix} \lambda_1 & \lambda_2\\\lambda_3&
                                                  \lambda_4\end{pmatrix}
  \begin{pmatrix} 0 & b\\ 0& d \end{pmatrix}=
  \begin{pmatrix} 0 & b\lambda_1+d\lambda_2
  \\  0  & b\lambda_3+d\lambda_4\end{pmatrix}\,,
 \ee
with the restriction
 \be
  |\det\Lambda|\,
  =\,|\lambda_1\lambda_4-\lambda_2\lambda_3|=1\,.
 \ee
Taking $\lambda_1=-d/{\rm gcd}(b;d)$ and $\lambda_2=b/{\rm
gcd}(b;d)$ we get
 \be
  \Lambda\,\g\,
  =\,\begin{pmatrix} -\frac{d}{{\rm gcd}(b;d)} & \frac{b}{{\rm gcd}(b;d)} \\
  \lambda_3&
                                                  \lambda_4\end{pmatrix}
  \begin{pmatrix} 0 & b\\ 0& d \end{pmatrix}=
  \begin{pmatrix} 0 & 0
  \\  0  & \frac{b\lambda_3+d\lambda_4}{{\rm gcd}(b;d)}\end{pmatrix}\,.
 \ee
Furthermore, since $\lambda_1=-d/{\rm gcd}(b;d)$ and
$\lambda_2=b/{\rm gcd}(b;d)$ are relatively prime, by the Bezout
identity there exist (non-unique)  $\lambda_3$ and $\lambda_4$ such
that
 \be\label{DetCon}
  b\lambda_3+d\lambda_4=1\,,
 \ee
providing
 \be
  \det\Lambda\,=\,\lambda_1\lambda_4-\lambda_2\lambda_3=1\,,
 \ee
so that $\Lambda\in GL_2(\IZ)$. As a consequence we obtain the coset
representative in the case $a=0$:
 \be
  \beta_{II}=\begin{pmatrix} 0 & 0\\0 &1\end{pmatrix}\,.
 \ee
The stabilizer of $\beta_{II}$ is found from the equation
$\Lambda\g=\g$, which reads
 \be
  \begin{pmatrix} \lambda_1 & \lambda_2\\\lambda_3&
                                                  \lambda_4\end{pmatrix}
  \begin{pmatrix} 0 & 0\\ 0& 1 \end{pmatrix}=
  \begin{pmatrix} 0 & \lambda_2
  \\ 0  & \lambda_4\end{pmatrix}\,.
 \ee
implying $\lambda_2=0$ and $\lambda_4=1$, which yields
 \be
  {\rm St}_{II}=\tilde{B}_-(\IZ)\,
  =\,\Big\{ \Big(\begin{smallmatrix}
  \pm 1&& 0\\&&\\m &&1\end{smallmatrix}\Big)\Big|\,m\in \IZ\Big\}\,.
 \ee

Second, we consider the case of $a\neq 0$ and $\lambda_3=0$. In this
case the condition $\det\g=ad=0$ implies $d=0$, then we get
 \be
  \Lambda\,\g\,
  =\,\begin{pmatrix} \lambda_1 & \lambda_2\\ 0&
                                                  \lambda_4\end{pmatrix}
  \begin{pmatrix} a & b\\ 0& 0 \end{pmatrix}=
  \begin{pmatrix} a\lambda_1 & b\lambda_1
  \\ 0  & 0 \end{pmatrix}\,.
 \ee
This allows to make $a>0$ and thus we arrive at  the following set of representatives
 \be\label{Coset2}
  \beta_{I}=\begin{pmatrix} a & b\\ 0& 0\end{pmatrix}, \qquad a\in
  \IZ_+\,\quad b\in \IZ\,.
 \ee
The stabilizer ${\rm St}_I$ of $\beta_I$ is found from the equation
 \be
  \begin{pmatrix} \lambda_1 & \lambda_2\\ \lambda_3&
                                                  \lambda_4\end{pmatrix}
  \begin{pmatrix} a & b\\ 0& 0 \end{pmatrix}=
  \begin{pmatrix} a\lambda_1 & b\lambda_1
  \\ a\lambda_3  & b\lambda_3 \end{pmatrix}\,,
 \ee
implying $\lambda_1=1$ and $\lambda_3=0$. Therefore we have
$|\lambda_4|=1$ providing the following
 \be
  {\rm St}_{I}=\tilde{B}_+(\IZ)\,=\,\Big\{ \Big(\begin{smallmatrix}
  \pm 1&& m\\&&\\0 &&1\end{smallmatrix}\Big)\Big|\,m\in \IZ\Big\}\,.
 \ee
This completes the proof. $\Box$

\begin{prop}
The following functional equation for the completed
  $GL_2$ zeta-function holds
 \be
  \xi^{GL_2}(1-s)\,=\,\xi^{GL_2}(s)\,,
 \ee
where the zeta-function is given by the analytic continuation  of
the integral \eqref{ALL123}:
 \be\label{ALL5}
  \xi^{GL_2}(s)\,=\int\limits_{GL_2(\IR)/GL_2(\IZ)}\!\!
  d\mu(g)\,\,|\det g|^{s}\,\,\Theta^{*}(0|\imath g^{\top}g)\,,
 \ee
and the truncated theta-function is given by
 \be\label{Theta25}
  \Theta^{*}(0|\imath g^{\top}g)\,
  =\sum_{\alpha\in\Mat^*(\IZ)}\!\!e^{-\pi \Tr(g^{\top}g\, \alpha\alpha^{\top})}\,.
  \ee
\end{prop}

\proof Note that theta-function \eqref{Theta25} possesses the
following symmetry
 \be
  \Theta^{*}(0|\imath \beta^{\top}g^{\top}g\beta)\,
  =\,\Theta^{*}(0|\imath gg^\top)\,,\qquad\beta\in GL_2(\IZ)\,,
 \ee
so the integration over the quotient space $GL_2(\IR)/GL_2(\IZ)$  is
well-defined. The theta-constant \eqref{Theta25} may be represented
as follows:
 \be
  \Theta^{*}(0|\imath g^{\top}g)\,
  =\,\Theta(0|\imath  g^{\top}g)\,-\,\Xi(0|\imath g^{\top}g)\,,
 \ee
where
 \be\label{Xi}
  \Xi(0|\imath g^{\top}g)\,
  =\!\!\sum_{\gamma\in\Mat_2^{(0)}(\IZ)}\!\!
  e^{-\pi\Tr(g^{\top}g\gamma\gamma^{\top})}\,,\\
  \Mat_2^{(0)}(\IZ)\,
  =\,\{\gamma\in\Mat(\IZ)\,|\quad\det \gamma=0\}\,.
 \ee

Notice that \eqref{Theta25} enjoys the following modular
transformation property \eqref{ModTr}:
 \be\label{ModTr5}
  \Theta(0|-T^{-1})=(\det T)^{1/2}\,\,\Theta(0|T)\,.
 \ee
Also note that for the involution $J:\,g\longmapsto g^{\tau}:=(g^{\top})^{-1}$, the following holds:
 \be
  g^{\top}g\longmapsto (g^{\top}g)^{-1}\,.
 \ee

Let us  split the integration domain in \eqref{ALL5} into two parts,
 \be
  GL_2(\IR)/GL_2(\IZ)\,=\,\CH_{\geq }\,\sqcup\,\CH_{\leq }\,,
 \ee
where
 \be
  \CH_{\geq }\,
  =\,\bigl\{g\in GL_2(\IR)/GL_2(\IZ)\,\bigl|\quad
  |\det g|\geq 1\bigr\}\,
  =\,\Mat^{\geq 1}_2(\IR)/GL_2(\IZ)\,, \\
  \CH_{\leq}\,
  =\,\bigl\{g\in GL_2(\IR)/GL_2(\IZ)\,\bigl|\,\,
  0<|\det g|\leq 1\bigr\}\,
  =\,\Mat^{\leq 1}_2(\IR)/GL_2(\IZ)\,,
 \ee
so that the  involution $J: g\mapsto g^\tau:=(g^{\top})^{-1}$
interchanges the two domains
 \be
  J\,:\quad\CH_{\geq}\longrightarrow \CH_{\leq}\,, \quad
  \CH_{\leq}\longrightarrow \CH_{\geq}\,.
 \ee
Also note that the transposition makes $J$ compatible with taking a 
quotient of $GL_2(\IR)$ over $GL_2(\IZ)$. Now we proceed with the
analog of the calculation for $GL_1$ for \eqref{ALL5}:
 \be
  \xi^{GL_2}(s)=\int\limits_{\CH_{\leq}}\!
  d\mu(g)\,\,|\det g|^{s}\,\,\Big(\Theta(0|\imath g^{\top}g)
  -\Xi(0|\imath g^{\top}g)\Big)\\
  + \int\limits_{\CH_{\geq}}\!
  d\mu(g)\,\,|\det g|^{s}\,\,\Big(\Theta(0|\imath g^{\top}g)
  -\Xi(0|\imath g^{\top}g)\Big)\,.
 \ee
Applying \eqref{ModTr5} to the first term of above expression we
obtain
 \be
  \xi^{GL_2}(s)\,
  =\,-\int\limits_{\CH_{\leq}}d\mu(g)\,\,|\det g|^{s}\,\,\Xi(0|\imath gg^\top)\\
  +\int\limits_{\CH_{\geq}}
  d\mu(g)\,\,|\det g|^{1-s}\,\,\Theta(0|\imath gg^\top)\\
  +\int\limits_{\CH_{\geq}}
  d\mu(g)\,\,|\det g|^{s}\,\,\Big(\Theta(0|\imath gg^\top)
  -\Xi(0|\imath gg^\top)\Big) \,,
 \ee
and then rewrite it as follows:
 \be
  \xi^{GL_2}(s)\,
  =\,-\int\limits_{\CH_{\leq}}d\mu(g)\,\,|\det g|^{s}\,\,\Xi(0|\imath g^{\top}g)\\
  +\int\limits_{\CH_{\geq}}d\mu(g)\,\,|\det g|^{1-s}\,\,\Xi(0|\imath g^{\top}g)\\
  +\int\limits_{\CH_{\geq}}d\mu(g)\,\,|\det g|^{1-s}\,\,\Big(\Theta(0|\imath g^{\top}g)-
  \Xi(0|\imath g^{\top}g)\Big) \\
  + \int\limits_{\CH_{\geq}}d\mu(g)\,\,|\det g|^{s}\,\,\Big(\Theta(0|\imath g^{\top}g)
  -\Xi(0|\imath g^{\top}g)\Big) \,.
 \ee
Next let us show that the sum of the first two terms in the
expression above,
 \be\label{Compare}
  -\int\limits_{\CH_{\leq}}d\mu(g)\,\,|\det g|^{s}\,\,\Xi(0|\imath g^{\top}g)\,
  +\int\limits_{\CH_{\geq}}d\mu(g)\,\,|\det g|^{1-s}\,\,\Xi(0|\imath
  g^{\top}g)\,,
 \ee
is symmetric with respect to $s\mapsto 1-s$. Indeed, by \eqref{Xi} the first term in \eqref{Compare} reads
 \be\label{term1}
  \int\limits_{\CH_{\leq}}\!d\mu(g)\,|\det g|^{s}\,\Xi(0|\imath g^{\top}g)\,
  =\int\limits_{\CH_{\leq}}\!d\mu(g)\,|\det g|^{s}\!\!\!\!
  \sum_{\gamma\in\Mat^{(0)}(\IZ)}\!\!\!e^{-\pi \Tr( g^{\top}g\gamma\gamma^{\top})}\,.
 \ee
Let us represent the integration domain in \eqref{term1} as follows
 \be\label{Hleq}
  \CH_{\leq}\times\Mat^{(0)}_2(\IZ)\,
  =\,\Big({\rm Mat}^{\leq 1}_2(\IR)/GL_2(\IZ)\Big)\times\Mat^{(0)}_2(\IZ)\\
  =\,\Big(\bigl(\Mat^{\leq 1}_2(\IR)/St_I\bigr)\times GL_2(\IZ)\,\beta_I\Big)
  \bigsqcup
  \Big(\bigl(\Mat^{\leq 1}_2(\IR)/St_{II}\bigr)\times GL_2(\IZ)\,\beta_{II}\Big)\,.
 \ee
Hence we re-write \eqref{term1} as follows:
 \be\label{term12}
  \int\limits_{\CH_{\leq}}\!d\mu(g)\,|\det g|^{s}\,\Xi(0|\imath g^{\top}g)\,
  =\,Z^{\leq 1}_I(s)\,+\,Z^{\leq 1}_{II}(s)\,. 
 \ee
Consider the second contribution
 \be\label{Z2}
  Z^{\leq 1}_{II}(s)\,
  =\!\!\!\int\limits_{\Mat^{\leq 1}_2(\IR)/J(St_{II})}\!\!\!\!
  d\mu(g)\,\,|\det g|^{s}\,
  \,e^{-\pi \Tr( g^{\top}g\beta_{II}\beta_{II}^{\top})}\,,
 \ee
then by \eqref{Coset0} we derive
 \be
  \beta_{II} \beta_{II}^{\top}\,=\,\begin{pmatrix} 0& 0\\0&1\end{pmatrix},
  \qquad \Tr(g^{\top}g\beta_{II} \beta_{II}^{\top})\,=\,(g^{\top}g)_{22}\,,
 \ee
and the stabilizer is given by \eqref{Stab}, 
 \be
  St_{II}= \tilde{B}_-(\IZ)=\Big\{
  \Big(\begin{smallmatrix}
  1&& 0\\&&\\m &&\pm1\end{smallmatrix}\Big)\Big|\,m\in \IZ\Big\}\,.
 \ee
Applying the Iwasawa decomposition (to the open part) of $\Mat^{\leq1}_2(\IR)$ gives (there
are lower-dimensional strata that we may neglect)
 \be
  \Mat^{\leq1}_2(\IR)\supset KAN_-\,, \quad g=kan_-\,,\quad
  d\mu(g)\,=\,dk\,dn_-\,a^{-2\rho}\,\frac{da_1\,da_2}{a_1a_2}\,,
 \ee
which implies
\be
g^{\top}g=n_-^{\top}a^2n_-=
\begin{pmatrix} 1 & n\\0& 1\end{pmatrix}
\begin{pmatrix} a_1^2 & 0\\0& a_2^2\end{pmatrix}
\begin{pmatrix} 1 & 0\\ n& 1\end{pmatrix}=\begin{pmatrix}
                                            a_1^2+n^2a_2^2 & na_2^2
                                            \\na_2^2& a_2^2\end{pmatrix}\,,
\ee
so that $(g^{\top}g)_{22}=a_2^2$. 
Therefore, considering  $N_-/St_{II}\simeq\IR/\IZ$ and $\int_Kdk=1$ , 
the integral \eqref{Z2} reduces to the following, for ${\rm Re}(s)>1$,  
 \be\label{Z20}
  Z^{\leq 1}_{II}(s)\,
  =\!\int\limits_{\IR/\IZ}\!dn_-\!\!
  \int\limits_{a_i\in \IR_+\atop{a_1a_2\leq 1}}\!\! 
  \frac{da_1\,da_2}{a_1a_2}\,\,|a_1a_2|^{s}\,a_1^{-2\rho_1} a_2^{-2\rho_2}
  e^{-\pi a_2^2}\\\
  =\int\limits_0^1da_2\,\,a_2^{s}\,\,e^{-\pi a_2^2}
  \int\limits_0^{a_2^{-1}} da_1\,\, a_1^{s-2}\,
  =\,\frac{1}{s-1}\int\limits_0^1da_2\,\,a_2\,\,e^{-\pi a_2^2}\,.
 \ee

Next let us calculate contribution of the other part of the contour \eqref{Hleq}  in \eqref{term12}:
 \be\label{Z1}
  Z^{\leq 1}_{I}(s)\,=\!\!
  \int\limits_{\Mat^{\leq 1}_2(\IR)/St_I}\!\! d\mu(g)\,\,|\det g|^{s}\,
  \,\sum_{\beta_I} e^{\imath \pi \Tr( g^{\top}g\beta_{I}\beta_{I}^{\top})}\,.
 \ee
where by \eqref{Coset3} we get, for $a\in\IZ_+,\,b\in\IZ$, 
 \be
  \beta_{I} \beta_{I}^{\top}=\begin{pmatrix} a^2+b^2& 0\\0&0\end{pmatrix},
  \qquad \Tr( gg^{\top}\beta_{I} \beta_{I}^{\top})=(g^{\top}g)_{11}(a^2+b^2)\,.
 \ee
and the stabilizer is given by \eqref{Stab}
 \be
  St_{I}= \tilde{B}_+(\IZ)\,
  =\,\Big\{ \Big(\begin{smallmatrix}
  \pm 1&& m\\&&\\0 &&1\end{smallmatrix}\Big)\Big|\,m\in \IZ\Big\}\,.
 \ee
Applying the Iwasawa decomposition of the open part of ${\rm Mat}_2(\IR)$
 \be
  N_-AK \subset {\rm Mat}_2(\IR), \quad g=kan_+\,,\quad
  d\mu(g)\,=\,dk\,dn_-\,a^{2\rho}\,\frac{da_1\,da_2}{a_1a_2}\,,
 \ee
(we might neglect lower-dimensional strata in the calculation of
the integral)
provides 
 \be
  g^{\top}g=n_+^{\top}a^2n_+
  =\begin{pmatrix} 1 & 0\\n& 1\end{pmatrix}
  \begin{pmatrix} a_1^2 & 0\\0& a_2^2\end{pmatrix}
  \begin{pmatrix} 1 & n\\ 0& 1\end{pmatrix}=\begin{pmatrix}
                                            a_1^2& na_1^2
                                            \\na_1^2& a_2^2+n^2a_1^2\end{pmatrix}\,,
 \ee
so that $(g^{\top}g)_{11}=a_1^2$\,.
Therefore, considering  $N_+/St_I\simeq\IR/\IZ$ and $\int_Kdk=1$ , 
the integral \eqref{Z1} boils down to the following, for ${\rm Re}(s)>1$:
 \be\label{Z10}
  Z^{\leq 1}_{I}(s)=\int\limits_{\IR/\IZ} dn
  \int\limits_{a_i\in \IR_+\atop{a_1a_2\leq 1}}\!\!\frac{da_1\,da_2}{a_1a_2}\,\,
  |a_1a_2|^s\,a_1^{2\rho_1} a_2^{2\rho_2}\sum_{a\in \IZ_+,b\in \IZ}\!\!e^{-\pi(a^2+b^2) a_1^2}\\
  =\int\limits_0^1da_1\,\,a_1^s\!\!\sum_{a\in \IZ_+,b\in \IZ}\!\!e^{-\pi(a^2+b^2) a_1^2}
  \int\limits_0^{a_1^{-1}}\!da_2\,\,a_2^{s-2}\\
  =\,\frac{1}{s-1}\int\limits_0^1da_1\,\,a_1\!\!\sum_{a\in \IZ_+,b\in \IZ}\!\!e^{-\pi(a^2+b^2) a_1^2}\,.
 \ee
Collecting together contributions \eqref{Z20} and \eqref{Z10} results in the following expression for the first term in \eqref{Compare}:
 \be
  -\int\limits_{\CH_{\leq}}d\mu(g)\,\,|\det g|^{s}\,\,\Xi(0|\imath g^{\top}g)\\
  =\,\frac{1}{1-s}\Big(\int\limits_0^1da_2\,\,a_2\,\,e^{-\pi a_2^2}\,
  +\int\limits_0^1da_1\,\,a_1\!\!\sum_{a\in \IZ_+,b\in \IZ}\!\!e^{-\pi(a^2+b^2) a_1^2}\Big)\,.
 \ee

Now we repeat the same calculations for the other term in \eqref{Compare}, 
 \be\label{term2}
  \int\limits_{\CH_{\geq}} d\mu(g)\,\,|\det g|^{\tilde{s}}\,\,\Xi(0|\imath
  g^{\top}g)\,,\qquad\tilde{s}:=1-s\,,\quad{\rm Re}(\tilde{s})<0\,,
 \ee
via splitting the integration domain similarly to \eqref{Hleq}: 
 \be\label{Hgeq}
  \CH_{\geq}\times\Mat^{(0)}_2(\IZ)\,
  =\,\Big({\rm Mat}^{\geq1}_2(\IR)/GL_2(\IZ)\Big)\times\Mat^{(0)}_2(\IZ)\\
  =\,\Big(\bigl(\Mat^{\geq1}_2(\IR)/St_I\bigr)\times GL_2(\IZ)\,\beta_I\Big)
  \bigsqcup
  \Big(\bigl(\Mat^{\geq1}_2(\IR)/St_{II}\bigr)\times GL_2(\IZ)\,\beta_{II}\Big)\,.
 \ee
All the computations are basically the same except for the 
integrals \eqref{Z20} and \eqref{Z10}.  Namely, contribution of the second
branch in \eqref{Hgeq} is given by, for ${\rm Re}(\tilde{s})<0$,
 \be\label{Z21}
  Z^{\geq1}_{II}(\tilde{s})\,
  =\!\int\limits_{\IR/\IZ}\!dn_-\!\!
  \int\limits_{a_i\in \IR_+\atop{a_1a_2\geq 1}}\!\! 
  \frac{da_1\,da_2}{a_1a_2}\,\,|a_1a_2|^{\tilde{s}}\,a_1^{-2\rho_1} a_2^{-2\rho_2}
  e^{-\pi a_2^2}\\\
  =\!\int\limits_0^1da_2\,\,a_2^{\tilde{s}}\,\,e^{-\pi a_2^2}
  \int\limits_{a_2^{-1}}^{\infty} da_1\,\, a_1^{\tilde{s}-2}\,
  =\,\frac{1}{1-\tilde{s}}\int\limits_0^1da_2\,\,a_2\,\,e^{-\pi a_2^2}\,.
 \ee
Similarly, contribution of the first
branch in \eqref{Hgeq} is given by, for ${\rm Re}(\tilde{s})<0$,
 \be\label{Z11}
  Z^{\geq 1}_{I}(\tilde{s})\,=\int\limits_{\IR/\IZ} dn
  \int\limits_{a_i\in \IR_+\atop{a_1a_2\geq 1}}\!\!\frac{da_1\,da_2}{a_1a_2}\,\,
  |a_1a_2|^{\tilde{s}}\,a_1^{2\rho_1} a_2^{2\rho_2}
  \sum_{a\in \IZ_+,b\in \IZ}\!\!e^{-\pi(a^2+b^2) a_1^2}\\
  =\!\int\limits_0^1da_1\,\,a_1^{\tilde{s}}\!\!
  \sum_{a\in \IZ_+,b\in \IZ}\!\!e^{-\pi(a^2+b^2) a_1^2}
  \int\limits_{a_1^{-1}}^{\infty}\!da_2\,\,a_2^{\tilde{s}-2}\\
  =\,\frac{1}{1-\tilde{s}}\int\limits_0^1da_1\,\,a_1\!\!
  \sum_{a\in \IZ_+,b\in \IZ}\!\!e^{-\pi(a^2+b^2) a_1^2}\,.
 \ee
Comparison of the  contributions of $\CH^{\geq 1}$ and $\CH^{\leq 1}$
 confirms the symmetry under $s\mapsto\tilde{s}:=1-s$ and therefore
 proves the  functional relation for the global/completed $GL_2$-zeta function $\xi^{GL_2}(s)$.

\vspace{5mm}


\noindent {\small {\bf A.A.G.} {\sl Laboratory for Quantum Field
Theory
and Information},\\
\hphantom{xxxx} {\sl Institute for Information
Transmission Problems, RAS, 127994, Moscow, Russia};\\
\hphantom{xxxx} {\it E-mail address}: {\tt anton.a.gerasimov@gmail.com}}\\
\noindent{\small {\bf D.R.L.} {\sl Laboratory for Quantum Field
Theory
and Information},\\
\hphantom{xxxx}  {\sl Institute for Information
Transmission Problems, RAS, 127994, Moscow, Russia};\\
\hphantom{xxxx} {\it E-mail address}: {\tt lebedev.dm@gmail.com}}\\
\noindent{\small {\bf S.V.O.} {\sl
 Beijing Institute of Mathematical Sciences and Applications\,,\\
\hphantom{xxxx} Huairou District, Beijing 101408, China};\\
\hphantom{xxxx} {\it E-mail address}: {\tt oblezin@gmail.com}}

\end{document}